\input amstex
\input epsf
\documentstyle{amsppt}
\magnification=1200
\raggedbottom

\widestnumber\key{DaKrSo}

\define \X{\frak X}

\define \E{\Cal E}
\define \a{\alpha}
\define \e{\varepsilon}

\define \ga{\gamma}
\define \de{\delta}
\define \De{\Delta}
\define \la{\lambda}
\define \La{\Lambda}
\define \ka{\varkappa}
\define \({\left(}
\define \){\right)}
\define \[{\left[}
\define \]{\right]}
\define \ov{\overline}
\define \un{\underline}

\define \n{\goth n}
\define \wt{\widetilde}
\define \Ker{\operatorname{Ker}}
\define \Tr{\operatorname{Tr}}
\define \supp{\operatorname{supp}}
\define \diam{\operatorname{diam}}
\define \pr{\prod_1^\infty\{0;1\}}
\define \pro{\prod_{-\infty}^\infty\{0;1\}}
\define \EE{\text{{\bf E}}}

\define\lr{\leftrightarrow}
\redefine \phi{\varphi}
\redefine \B{\frak B}
\redefine \S{\Cal S}

\topmatter
\title  Ergodic properties of Erd\"os measure,
the entropy of the goldenshift, and related problems
\endtitle

\date December 2, 1997 \enddate

\author Nikita SIDOROV and Anatoly VERSHIK \endauthor
\address St.~Petersburg Department of Steklov Mathematical Institute,
Fontanka 27, St.~Petersburg 191011, Russia \endaddress
\email sidorov$\@$pdmi.ras.ru,
vershik$\@$pdmi.ras.ru \endemail
\rightheadtext{Erd\"os measure and the goldenshift}

\abstract
We define a two-sided analog of Erd\"os measure on the
space of two-sided expansions with respect to
the powers of the golden ratio, or, equivalently,
the Erd\"os measure on the 2-torus. We construct
the transformation (goldenshift) preserving both Erd\"os
and Lebesgue measures on $\Bbb T^2$ which is the induced
automorphism with respect to the ordinary shift (or the corresponding
Fibonacci toral automorphism) and proves to be Bernoulli
with respect to both measures in question. This provides a direct
way to obtain formulas for the entropy dimension of the
Erd\"os measure on the interval, its entropy in the sense
of Garsia-Alexander-Zagier and some other results.
Besides, we study central measures on the Fibonacci
graph, the dynamics of expansions and related questions.
\endabstract
\dedicatory  To the memory of Paul Erd\"os \enddedicatory
\endtopmatter

\document

\footnote" "{Supported in part
by the INTAS grant 93-0570 and RFBR grant 96-0100676.
The first author was supported by the
French foundation PRO MATHEMATICA. The first
author expresses his gratitude to l'Institut
de Math\'ematiques de Luminy for support during
his stay in Marseille in 1996-97. The second author is grateful
to the University of Stony Brook for support during his visit
in February-March 1996 and to the Institute for Advanced studies
of Hebrew University for support during his being there in 1997.
The previous version of this paper appeared as the Stony Brook
preprint IMS 96-14.}

\head 0. Introduction \endhead

Among numerous connections between ergodic theory and the metric
theory of numbers, the questions related to algebraic irrationals,
expansions associated with them and ergodic properties
of related dynamical systems, are of special interest.
The simplest case, i.e. the golden ratio, the Fibonacci
automorphism etc., has served as a deep source of problems and
conjectures.

In 1939 P.~Erd\"os \cite{Er} proved in particular
the singularity of the measure
on the segment which is defined as the one corresponding to
the distribution of the random variable $\sum_1^\infty \e_k\la^{-k}$,
with $\la$ being the larger golden ratio
and $\e_k$ independently taking the values 0 and 1 (or $\pm1$)
with probabilities $\frac12$ each. We think it is natural to
call this measure the {\it Erd\"os measure}.
This work gave rise to many
publications and numerous generalizations (see, e.g., \cite{AlZa}
and references therein). Nevertheless, little attention was paid
to the dynamical properties of this natural measure.
{\it The aim of this paper is to begin studying dynamical properties
of Erd\"os measure and its two-sided extension}. We
\roster
\item define the two-sided generalization of Erd\"os measure (Section~1);
\item introduce a special automorphism (``goldenshift'') which preserves
Erd\"os measure
and which is a Bernoulli automorphism with a natural generator
with respect to Erd\"os and Lebesgue measures (Section~2);
\item compute the entropy of this automorphism and prove that it is
related to the pointwise dimension of the Erd\"os
measure defined in \cite{Y}
and Garsia's entropy considered in \cite{AlZa} (Section~3);
\item discuss the connections with some properties of the Fibonacci
graph, its central measures and the adic transformation on it (see
Appendix~A);
\item define a new kind of expansions corresponding to the
goldenshift (see Appendix~B).
\endroster
We will describe all this in more detail below.

Several years ago certain connections between symbolic dynamics of toral
automorphisms and arithmetic expansions associated with their
eigenvalues were established. The first step in this direction
was also related to the golden ratio (see \cite{Ver5}) and led to
a natural description of the Markov partition in terms of the arithmetic
of the 2-torus and homoclinic points of the Fibonacci automorphism.
The main idea was to consider the natural extension
of the shift in the sense of ergodic theory and the adic transformation
on the space of one-sided arithmetic expansions and to identify
the set of {\it two-sided expansions}
with the 2-torus (see also \cite{Ver6}, \cite{Ver3}, \cite{KenVer}).

In the present paper we use the same idea for a detailed study
of the Erd\"os measure. Namely, we define the two-sided Erd\"os
measure as a measure on the space of expansions infinite to
both sides and
identify it with a measure on the 2-torus; in the same way,
Lebesgue measure on the 2-torus can be considered as a two-sided
version of the Markov invariant measure on the corresponding
Markov compactum -- see below and \cite{Ver5}, \cite{Ver6}.
We study the properties
of the ordinary shift and the {\it goldenshift}
as a transformation on the space of expansions introduced
by means of the notion of {\it block}.
The goldenshift turns out to preserve both Lebesgue and Erd\"os measures,
both being Bernoulli in the natural sense with respect to
the goldenshift; this is one of the main results of the paper
(Theorem~2.7).
By the way, this immediately yields a proof of the Erd\"os theorem
on the singularity of Erd\"os measure.
Moreover, the two-sided goldenshift is an induced automorphism
for the Fibonacci automorphism of the torus.
Other important consequences
of our approach follow from the fact that the entropy of the goldenshift
is directly related to the entropy of Erd\"os measure in the
sense of Garsia and Alexander-Zagier, i.e. to the entropy of the
random walk with equal transition measures on the Fibonacci
graph (Theorem~3.3).

In \cite{AlZa} it was attempted to compute the
entropy of Erd\"os measure as the infinite convolution of
discrete measures, which was in fact introduced by A.~Garsia
\cite{Ga} in a more general situation. Note that it proved to
be the entropy of a random walk on the Fibonacci graph. In
\cite{LePo} the authors compute the dimension of the
Erd\"os measure on the interval in the sense of L.-S.~Young \cite{Y} and
relate a certain two-dimensional dynamics to it.

 Finally, making use of a version of Shannon's theorem
for random walks (see \cite{KaVer}) yields the value of the dimension of
the Erd\"os measure in the sense of Young (Theorem~3.7).

Thus, the dynamical viewpoint for arithmetic expansions and for
measures related to them provides new information and an essential
simplification of computations of the invariants invloved.
One may expect that the methods of this paper will apply
to more general algebraic irrationals and also
to some nonstationary problems.

The contents of the present paper are as follows. In Section~1
we present auxiliary notions (canonical expansions and others)
and give the main definitions (Erd\"os measure on the interval
and the 2-torus, normalization, the Markov measure
corresponding to Lebesgue measure etc.). In addition, we deduce
some preliminary facts about the one-sided and two-sided Erd\"os
measures. In Section~2 we study the combinatorics of
blocks in terms of canonical expansions, introduce the notion
of the goldenshift (both one-sided and two-sided) and prove
its Bernoulli property with respect to Lebesgue and Erd\"os
measures. Section~3 contains our main results on the entropy
and dimension of the Erd\"os measure and their relationships
to the random walk on the Fibonacci graph.

In the appendices we consider some related problems. Namely,
in Appendix~A the combinatorial and algebraic theory of the
Fibonacci graph is presented. In particular, we describe the
ergodic central measures on this graph and the action
of the {\it adic} transformation, which is defined as
the transfer to the immediate successor in the sense of
the natural lexicographic order (in our case it is just
 the next expansion of a given real in the sense of
the natural ordering of the expansions). We study the metric
type of the adic transformation with respect to the ergodic
central measures. In Appendix~B we consider arithmetic block
expansions of almost all points of the interval. The interest
in them is due to the fact that the ``digits'' of the block
expansions are independent with respect both to Erd\"os
and Lebesgue measures.
Note that there are some pecularities caused by the difference
between the one-sided and two-sided shifts. For instance,
the one-sided Erd\"os measure is only quasi-invariant under
the one-sided shift, while the two-sided measure is shift-invariant.
In Appendix~C the densities of the Erd\"os measure
with respect to the shift and to the rotation by the golden
ratio are computed by means of blocks.
Finally, in Appendix~D another proof of Alexander-Zagier's formula
for the entropy is given. It is worthwhile because of its
connection with the geometry of the Fibonacci graph.

The authors express their gratitude to B.~Solomyak for the fruiful
discussions and to the referees for their helpful remarks and
suggestions.

\head 1. Erd\"os measure on the interval and on the 2-torus. \endhead

\subhead 1.1. Canonical expansions \endsubhead
Let
$\Sigma=\pr$ endowed with the one-sided Bernoulli
shift $\sigma$ and let
$X\subset \Sigma$ be the stationary Markov compactum with
the transition matrix
$\(\smallmatrix 1&1\\1&0\endsmallmatrix\)$, i.e. the set
$X=\{(\e_1\e_2\dots)\in\Sigma:
\e_{i}\e_{i+1}=0,\ i\ge 1\}$ endowed with the topology of pointwise
convergence. Let next $\la=\frac{\sqrt5+1}2$ and
$L:X\to [0,1]$ be the mapping acting by the formula
$$
L(\e_1\e_2\dots):=\sum_{k=1}^\infty \e_k\la^{-k}.\tag 1.1
$$
It is well known that $L$ is one-to-one, except for a countable number
of sequences whose tail is of the form
$0^\infty$ or $(01)^\infty$.
The inverse mapping $L^{-1}$ is specified with the help
of the {\it greedy algorithm}. Namely, let $Tx=\{\la x\}$ and
$$
\e_k=[\la T^{k-1}x], \quad k\ge1.
$$

We call the constructed sequence $(\e_1(x) \e_2(x)\dots)$ the {\it
canonical expansion} of $x$.
Note that usually the canonical expansions are called
{\it $\beta$-expansions} (for $\beta=\la$). They were
introduced in \cite{Re} and \cite{Ge} and thoroughly
studied in \cite{Pa}.

Note that the mapping $L$ can be naturally extended to $\Sigma$, and
let $L'$ stand for this extension. However,
$L'(\Sigma)=[0,\la]$, that is why we define the projection
$\pi:\Sigma\to[0,1]$ by the formula
$\pi(x):=\la^{-1}L'(x)$. It will be frequently
used below. Note that $\pi$ is not one-to-one.

\subhead 1.2. The Markov measure on $X$ \endsubhead
The transformation $T:(0,1)\to(0,1)$ is transferred by $L$ to the Markov
compactum $X$ and acts as the one-sided shift $\tau$:
$$
\tau(\e_1\e_2\e_3\dots):=\e_2\e_3\dots
$$

Thus, we have $\tau=L^{-1}TL$.
The transformation $T$ has been thoroughly studied, and
it was shown that there exists the $T$-invariant
measure $m'$ equivalent to
Lebesgue measure $m_1$. Its density is given by the formula
$$
\rho(x)=\frac{dm'}{dm_1}=\cases \la^2/\sqrt5, & 0<x\le \la^{-1} \\
               \la/\sqrt5, & \la^{-1}<x\le 1
        \endcases
$$
(see, e.g, \cite{Ge}, \cite{Pa}).
The corresponding Markov measure $L^{-1}m'$ on $X$ is the one
with the stationary initial
distribution $\(\smallmatrix \la/\sqrt5 \\ \la^{-1}/\sqrt5 \endsmallmatrix\)$
and the transition probability matrix $\(\smallmatrix \la^{-1} &
\la^{-2} \\ 1&0\endsmallmatrix \)$.
The $L$-preimage of the Lebesgue measure $m$ on $X$
differs from this stationary Markov measure only by its initial
distribution $\(\smallmatrix \la^{-1} \\ \la^{-2}\endsmallmatrix\)$.
Note that for the adic transformation on $X$ (for the definition
see \cite{Ver2}
or Appendix~A) with the alternating ordering on the paths
the latter measure is unique invariant, as this adic transformation turns
into the rotation by the angle $\la^{-1}$ under the mapping
$L$ (for more details see \cite{VerSi}).

\subhead 1.3. Erd\"os measure and normalization \endsubhead
Let us define the {\it Erd\"os}
measure.
By definition, the continuous Erd\"os measure $\mu$ on the unit interval is
the infinite convolution $\vartheta_1*\vartheta_2*\dots$, where
$\supp\vartheta_n=\{0,\la^{-n-1}\}$, and $\vartheta_n(0)=
\vartheta_n(\la^{-n-1})=\frac12$
(see \cite{Er}).

We are going to specify this measure more explicitly.
Let $p$ denote the product measure with the equal
multipliers $\(\frac12, \frac12\)$ on the compactum $\Sigma$.
Then it is easy to see that $\mu=\pi(p)$.

We are also interested in the specification of the Erd\"os measure
on the Markov compactum $X$. Of course, it is just $L^{-1}\mu$;
however, it is worthwhile to introduce a direct mapping.
\definition{Definition} Let $x\in\Sigma,\ x=\{x_k\}_{k=1}^\infty$;
we define $[0,1]\ni c(x)=\sum_{k=1}^\infty x_k\la^{-k-1}=
\sum_{k=1}^\infty \e_k\la^{-k}$, where $\{\e_k\}$ is the canonical
expansion of $c(x)$. We define
$$
\n(x):=\e=\{\e_k\}_{k=1}^\infty.
$$
The mapping $\n:\Sigma\to X$ is called {\it normalization}.
\enddefinition
We will also describe the mapping $\n$ directly avoiding the expansion
of a number from $[0,1]$. Namely, let $x=(x_1,x_2,\dots)\in\Sigma$;
we put $x_0=0$
and look for the first
occurrence of the triple 011, after which we
replace it by 100. The next step is the same, i.e. we
return to the zero coordinate and start from there until we meet
again 011, etc. It is easy to see that the process leads to stabilization
at a normalized sequence. Note that this algorithm is rather rough,
as it is known that there exists a finite automaton
carrying out the process of normalization faster (see, e.g., \cite{Fr}).
It is obvious that this definition is equivalent to the one given
above.

Now we can define Erd\"os measure also on $X$ as the image of the
product measure $p=\prod_1^\infty\{\frac12,\frac12\}$:
$$
\mu=\n(p)
$$
(we preserve the same notation for $X$ as for $[0,1]$).

Below we will see that this definition of the Erd\"os measure is
not very suitable for deducing its dynamical properties
(the quasi-invariance under $T$ and the rotation by the golden ratio,
etc.); we will give another definition related to the two-sided
theory.

We begin with its self-similar property which is typical for
this measure and completely characterizes it.

\proclaim{Lemma 1.1} The Erd\"os measure $\mu$ on the interval~$[0,1]$
satisfies the following self-similar relation:
$$
\mu E=\cases \frac12\mu(\la E), & E\subset [0, \la^{-2}) \\
             \frac12 (\mu(\la E) + \mu(\la E-\la^{-1})), &
                                            E\subset [\la^{-2}, \la^{-1}) \\
             \frac12 \mu(\la E-\la^{-1}), & E\subset [\la^{-1}, 1]
       \endcases
$$
for any Borel set $E$.
\endproclaim
\demo{Proof} Let $F_1=1, F_2=2,\dots$ be the sequence of Fibonacci numbers.
Let $f_n(k)$ denote the number
of representations of a nonnegative integer $k$ as a sum of
not more than $n$ first Fibonacci numbers. We first show that for $n\ge 3$,
$$
f_n(k)= \cases f_{n-1}(k), & 0\le k\le F_n-1 \\
               f_{n-1}(k) + f_{n-1}(k-F_n), & F_n\le k\le F_{n+1}-2 \\
               f_{n-1}(k-F_n), & F_{n+1}-1 \le k \le F_{n+2}-2.
        \endcases  \tag 1.2
$$
To prove this, we represent $f_n(k)$ as $f_n(k)=f_n'(k)
+f_n''(k)$ for each $k<F_{n+2}-1$, where $f_n'(k)$
is the number of representations with  $\e_n=0$,
and $f_n''(k)$ is the number with $\e_n=1$. Obviously, if $k\le F_n-1$, then
$k=\sum_1^n \e_jF_j = \sum_1^{n-1} \e_jF_j$, whence $f_n(k)=
f_n'(k)$. If $F_{n+1}-1 \le k\le F_{n+2}-2$, then
$f_n(k)=f_n''(k)$. In the case $F_n\le k\le F_{n+1}-2$, obviously,
$f_n'(k)>0,\, f_n''(k)>0$. It remains to note only that
$f_n'(k)=f_{n-1}(k)$, and $f_n''(k)=f_{n-1}(k-F_n)$.

Now from (1.2),
and from the definition of the Erd\"os measure it follows
that  if, say, an interval $E\subset [0, \la^{-2})$, then
$$
\mu E=\lim_{n\to\infty} \sum_{k:\frac k{F_{n+2}} \in E} \frac {f_n(k)}{2^n}=
\frac12\lim_{n\to\infty} \sum_{k:\frac k{F_{n+1}} \in \la E}
\frac {f_{n-1}(k)}{2^{n-1}}=\frac12 \mu(\la E).
$$
The other cases are studied in the same way.\qed
\enddemo
\remark{Remark} The Erd\"os measure $\mu$ (as a Borel measure)
 is completely determined
by the above self-similar relation. Indeed, by induction one can determine
its values for any interval~$(a,b)$ with $a,b\in(\Bbb Z+\la\Bbb Z)\cap [0,1]$.
\endremark
\proclaim{Corollary 1.2} $\mu(0,\la^{-2})=\mu(\la^{-2},\la^{-1})=
\mu(\la^{-1},1)=\frac13$.
\endproclaim
The next step consists in introducing a two-sided analog of the Erd\"os
measure, which will lead to the one-sided shift-invariant measure
equivalent to $\mu$.

\subhead 1.4. Two-sided theory \endsubhead
Consider the two-sided space $\wt\Sigma=\pro$ and
its subset the two-sided Markov compactum
$\wt X=\{\{\e_k\}_{-\infty}^\infty:
\e_k\in\{0,1\},\ \e_k\e_{k+1}=0,\ k\in\Bbb Z\}$.\footnote{Henceforth
the sign \it tilde \rm will always stand for the two-sided objects.}
To define the two-sided Erd\"os measure on $\wt X$, we are going to
construct a two-sided analog of normalization. Furthermore,
we will use an arithmetic mapping from $\wt X$ onto $\Bbb T^2$
which semiconjugates the two-sided Markov shift and the Fibonacci automorphism
of the torus in order
to specify Erd\"os measure on the 2-torus and to
study its properties (see item~1.6).

Let $\wt\sigma$ denote the two-sided shift on $\wt\Sigma$, i.e.
$(\wt\sigma x)_k=x_{k+1}$, and let $\wt\tau$ stand for the two-sided
shift on the Markov compactum $\wt X$. We denote by $\wt m$ the stationary
two-sided Markov
measure on the compactum $\wt X$ with the invariant initial distribution
$\(\smallmatrix \la/\sqrt5 \\ \la^{-1}/\sqrt5 \endsmallmatrix\)$
and the transition probability matrix $\(\smallmatrix \la^{-1} &
\la^{-2} \\ 1&0\endsmallmatrix \)$. As is well known, $\wt m$
is the unique measure of maximal entropy for the shift $\wt\tau$.

There is an important action of $\Bbb Z^2=\Bbb Z+\Bbb Z$ on $\wt X$.
Namely, let $w_0,\ w_1$ be the generators, i.e. $\Bbb Z^2=\{nw_0+mw_1\mid
n,m\in\Bbb Z\}$. Let us describe the action $A_g:\wt X\to\wt X,\
g\in\Bbb Z^2$. First, $A_{w_1}=\wt\tau^{-1} A_{w_0}\wt\tau$,
and $A_{w_0}$ is addition by 1 in the sense of the arithmetic of $\wt X$.
More precisely, considering $\wt X$ as the set of formal series
$\{\sum_{-\infty}^\infty \e_nw_n\mid \{\e_n\}\in\wt X\}$,
we define $w_n+w_{n+1}=w_{n-1}$, and $2w_n=w_{n-1}+w_{n+2}$
(implying the representation $w_n\leftrightarrow \la^{-n}$),
whence the operation $\{\e_n\}\mapsto\{\e_n\}+w_0$ is well
defined for $\wt m$-a.e. $\{\e_n\}\in\wt X$, as well as
the sum of a.e. pair of sequences (see \cite{Ver5},
\cite{Ver3}).

\proclaim{Proposition} {\rm (see \cite{Ver5}, \cite{Ver3}).}
The measure $\wt m$ is the unique Borel measure invariant
under the action of $\Bbb Z^2$ decribed above.
\endproclaim
\remark{Remark} We define the following
identification of some pairs of points in $\wt X$ whose measure
$\wt m\times\wt m$
is 0 in order to turn $\wt X$ into the additive group.
Let the equivalence relation $\sim$ be defined as follows:
$$
\align
(*1000000\dots)&\sim (*010101\dots), \\
(\dots 01010100*)&\sim (\dots 10101001*),
\endalign
$$
where $*$ denotes an arbitrary (but the same for both sequences)
tail starting at the same term. Besides,
extending the equivalence relation~$\sim$
by continuity, we get $(0.1)^\infty\sim(1.0)^\infty\sim
0^\infty$
(henceforth the point will mark the
border between the negative and nonnegative coordinates of a sequence).
Now one can easily check that the set $\wt X'=\wt X/\sim$ is a compact
connected group in addition (see \cite{Ver5}, \cite{Ver3}).

\example{Example} Here is an example of subtraction in $\wt X$:
$-\la=\la^2+\la^4+\la^6+\ldots=1+\la^3+\la^5+\la^7+\ldots$, both
sequences representing one and the same element of $\wt X'$. Similarly,
for any sequence $\e\in\wt X$ finite to both sides, $-\e$ is
a pair of sequences finite to the right and {\it cofinite} to
the left, i.e. with the left tail $(01)^\infty$.
\endexample
\endremark

The purpose for the definition of
the operation of addition will be explained in item~1.6,
where the automorphism $(\wt X,\wt m,\wt\tau)$ will be related to the 2-torus.

Following the one-sided framework, we are going to define
the two-sided generalization of the operation of one-sided normalization.
Namely, we define the {\it two-sided normalization} $\wt\n$
as the mapping from $\wt\Sigma$ to $\wt X$.
Consider the set
$\wt\Sigma'\subset \wt\Sigma$ defined as follows:
$$
\wt\Sigma' = \left\{x\in\wt\Sigma: \#\{k<0:x_k=0\}=\infty\right\}.
$$
On the set $\wt\Sigma'$ we will define two-sided normalization.
Let $x\in\wt\Sigma'$ and $0\ge k_1>k_2>\dots,\ x_{k_i}=0,\ x_k\neq0,
k\neq k_i$ for all $i$. We set $x^{(0)}=x$. Let
$
\n\(\{x_i^{(0)}\}_{k_1+1}^\infty\)=\{\e_i^{(1)}\}_{i=k_1}^\infty.
$
We define
$$
x_i^{(1)}=\cases \e_i^{(1)} & i\ge k_1 \\
                 x_i        & i<k_1.
          \endcases
$$
By induction, let
$
\n\(\{x_i^{(n-1)}\}_{k_n+1}^\infty\)=\{\e_i^{(n)}\}_{i=k_n}^\infty.
$
Then, by definition,
$$
x_i^{(n)}=\cases \e_i^{(n)} & i\ge k_n \\
                 x_i        & i<k_n.
          \endcases
$$
Obviously, the process leads to the stabilization of $x_i^{(n)}$
in $n$.
\definition{Definition} The two-sided normalization $\wt\n$ at
$x\in\wt\Sigma'$ is defined as follows:
$$
\wt\n(\{x_i\}_{-\infty}^\infty)=\lim_{n\to\infty}\{x_i^{(n)}\}.
$$
\enddefinition

\definition{Definition} The two-sided Erd\"os measure $\wt\nu$ on the
Markov compactum $\wt X$ is
the image under the mapping $\wt\n$ of the
measure $\wt p$, which is the product of infinite factors
$\(\frac12, \frac12\)$ on the full compactum
$\wt\Sigma$.
\enddefinition
Since the set $\wt\Sigma'$ has full measure $\wt p$, we have
the homomorphism of the measure spaces:
$$
\wt\n:(\wt\Sigma, \wt p)\to(\wt X,\wt\nu).
$$
Let now
$$
\align
\rho&: \wt\Sigma\to\Sigma,\ \rho(\{x_n\})=(x_1, x_2,\dots),\\
\rho'&:\wt X\to X,\ \rho(\{\e_n\})=(\e_0, \e_1,\dots)
\endalign
$$
be the projections.

Let us write the following diagram:
$$
\matrix
\format\r&\,\,\,\c&\,\,\r&\,\,\,\c&\,\,\r&\,\,\,\c&\,\,\r\\
\Sigma & @<{\quad\sigma\quad} << &
\Sigma & @<{\quad\rho\quad}   << &
\wt\Sigma &@<{\quad\wt\sigma\quad}<< &
\wt\Sigma \\
 & \ssize\text{not commutes} && \ssize\text{not commutes}&&
 \ssize\text{commutes} & \\
{\ssize\n}\!\Big\downarrow & &{\ssize\n}
\!\Big\downarrow & &
{\ssize\wt\n}\!\Big\downarrow & &
{\ssize\wt\n}\!\Big\downarrow \\
X & @<{\quad\tau\quad}<< & X & @<{\quad\rho'\quad}<< &
\wt X & @<{\quad\wt\tau\quad}<< & \wt X
\endmatrix
$$

Note that
$$
\rho'\wt\n \neq \n\rho,\,\,\tau\n\neq\n\sigma.
$$
This is the reason why the one-sided theory has some difficulties.
The two-sided theory is more natural for this purpose
(the right part of the diagram does commute).

\proclaim{Propostion 1.3}
The two-sided Erd\"os measure $\wt\nu$ is invariant under
the two-sided shift, i.e. $\wt\tau\wt\nu=\wt\nu$
(cf. the one-sided case, where it does not take place).
\endproclaim
\demo{Proof} From the above specification of the mapping $\wt\n$
it follows that
$$
\wt\n\wt\sigma=\wt\tau\wt\n,\tag 1.3
$$
hence $\wt\nu(\wt\tau^{-1}E)=\wt p(\wt\n^{-1}\wt\tau^{-1}E)=
\wt p(\wt\sigma^{-1}\wt\n^{-1}E)=(\wt\sigma\wt p)(\wt\n^{-1}E)=
\wt p(\wt\n^{-1}E)=\wt\nu(E)$ for any Borel set $E\subset \wt X$.
\enddemo

\proclaim{Proposition 1.4} For any cylinder $\wt C=(\e_1=i_1,\dots,\e_r=i_r)
\subset \wt X$, its measure $\wt\nu$ is strictly positive.
\endproclaim
\demo{Proof} It follows from the direct specification of the two-sided
normalization described above that for the cylinder
$C'=(\e_0=0,\e_1=i_1,\dots, \e_r=i_r,\e_{r+1}=0,\e_{r+2}=0)\subset\wt\Sigma$,
we have $\wt\n^{-1}(\wt C)\supset C'$,
whence, by definition of the Erd\"os measure,
$\wt\nu(\wt C)\ge 2^{-r-3}$.
\enddemo

\proclaim{Proposition 1.5} The two-sided shift $\wt\tau$ on the Markov
compactum $\wt X$ with the two-sided Erd\"os measure
is a Bernoulli automorphism.
\endproclaim
\demo{Proof} We observe that this dynamical system is a factor
of the Bernoulli shift $\wt\sigma:\wt\Sigma\allowmathbreak\to\wt\Sigma$
with the product
measure $\(\frac12, \frac12\)$ (see relation~(1.3))
and apply the theorem due to D.~Ornstein
\cite{Or} on the Bernoullicity of all Bernoulli factors.
\enddemo
\remark{Remark} Note that the measure $\wt\nu$
is not Markov on the compactum $\wt X$.
It would be interesting to prove that the two-sided
Erd\"os measure is a Gibbs
measure for a certain natural potential.
\endremark

\subhead 1.5. Dynamical properties of Erd\"os measure
\endsubhead
Let the measure $\nu$ on the one-sided compactum $X$ be defined as
the projection of the two-sided Erd\"os measure $\wt\nu$. In other
words, the dynamical system $(\wt X, \wt\nu, \wt\tau)$ is
the natural extension of $(X, \nu, \tau)$.
We recall that it means by definition that for any cylinder
$C=(\e_1=i_1,\dots,\e_k=i_k)\subset X$ its measure $\nu$ equals
$\wt\nu(\wt C)$, where $\wt C=(\e_1=i_1,\dots,\e_k=i_k)\subset \wt X$.
Below we will denote the $L$-image of $\nu$ on the interval $[0,1]$
by the same letter.

\proclaim{Proposition 1.6}
The measure $\nu$ is $\tau$-invariant and ergodic.
\endproclaim
\demo{Proof}
The $\tau$-invariance of $\nu$ is a consequence of the $\wt\tau$-invariance
of $\wt\nu$. Furthermore, since the automorphism $(\wt X,\wt\nu,\wt\tau)$
is Bernoulli, it is ergodic, thus, the endomorphism $(X,\nu,\tau)$
is also ergodic.
\enddemo

We are going to establish a relation between the mappings
$\rho'\wt\n$ and $\n\rho$ in order to prove the equivalence
of the measures $\mu$ and $\nu$ on the interval.
Note first that in these terms
$\mu=(\n\rho)(\wt p)$, and $\nu=(\rho'\wt\n)(\wt p)$.
We also note that the fact that $\mu\prec\nu$ is shown
easily, while in the opposite direction it is not straightforward.
The following claims yield the same proofs for both sides.

\proclaim{Lemma 1.7}
There exists a subset of $\wt\Sigma$ of full measure $\wt p$ and its
countable partition into the
sets $\{E_k\}_{k=1}^\infty$ and the corresponding set $\{f_k\}_{k=1}^\infty$
of finite sequences in $\wt\Sigma$ such that
$$
(\rho'\wt\n)(x)=(\n\rho)(x+f_k),\quad x\in E_k\tag 1.4
$$
with group addition in the set $\wt\Sigma=\prod_{-\infty}^\infty \Bbb Z/2$.
Similarly,
there exists a $\wt p$-a.e. partition of $\,\wt\Sigma$ into
the
sets $\{D_k\}_{k=1}^\infty$ and the corresponding set $\{g_k\}_{k=1}^\infty$
of finite sequences in $\wt\Sigma$ such that
$$
(\n\rho)(x)=(\rho'\wt\n)(x+g_k),\quad x\in D_k.\tag 1.5
$$
\endproclaim
\remark{Remark} The actions $x\mapsto x+f_k,\ x\mapsto x+g_k$ change
finitely many coordinates of $x$. We need only this property.
\endremark
\demo{Proof} Both assertions are proved in the same
way. Let us prove the first one.
The idea of the proof is based on the fact that a $\wt p$-typical
sequence from $\wt\Sigma$ can be splitted into finite pieces so that
its normalization splits into the concatenation of the corresponding
normalizations.

We assume that $x$ has two successive zero coordinates
with negative indices and four successive zero coordinates with
positive indices. So, let $x=(A_0 00 A_1 0000 A_2)$, where $A_0$
and $A_2$ are infinite and $A_1$ is a finite fragment of $x$ containing
$x_0$ and $x_1$. Then
$\wt\n(x)$ is the concatenation of the normalizations of its pieces
$A_00, 0A_1000$ and $0A_2$,
whence
$$
(\rho'\wt\n)(x)=(C0\n(0A_2)) \tag1.6
$$
with a certain finite admissible word $C$ ending with two zeroes.
We have a countable number of possibilities for $C$.
Let $E_k:=E_C$ be defined as the set of $x\in\wt\Sigma'$ such that
relation~(1.6) holds with some $A_2$.

To construct $f_k$, we
consider two cases. If $C$ begins with 0, i.e. if $C=0C'$, then
we set $x':=x+f_k=(A_00\dots0\mid C'00A_2)$, where ``$\mid$'' denotes
the border between positive and nonpositive coordinates.
For such an $x'$ relation~(1.4) is satisfied. If $C$ begins
with 1, i.e. if $C=(10)^j 0C'$ for some $j\ge1$ and admissible
$C'$, then we set $x':=(A_00\dots0\mid 10(11)^{j-1}C'00A_2)$.
The proof is complete.
\enddemo

Relations (1.4) and (1.5) together yield the main assertion.

\proclaim{Proposition 1.8} The measures $\mu$ and $\nu$ are equivalent.
\endproclaim
\demo{Proof} Consider the group which acts on $\wt\Sigma$
by adding the finite sequences ($=$ by changing a finite
number of coordinates). Since the action of this group
preserves the measure $\wt p$,
 we obtain from relation~(1.4)
$(\rho'\wt\n)(\wt p) \prec (\n\rho)(\wt p)$, and from
relation~(1.5) we get $(\n\rho)(\wt p)\prec(\rho'\wt\n)(\wt p)$,
whence $\nu\prec\mu$, and $\mu\prec\nu$.\qed
\enddemo

\remark{Remark} It is possible to show that there exist
two positive constants $C_1$ and $C_2$ such that
$$
C_1\mu(E) \le \nu(E) \le -C_2\mu(E)\log\mu(E)
$$
for any Borel $E$.
Note that the right estimate is attainable, for instance,
at the sequence of sets $E=E_n=(0,\la^{-n})$, as by Lemma~1.1
and Corollary~1.2,
$\mu(0,\la^{-n})=\frac43\cdot 2^{-n}$, while
$\nu(0,\la^{-n})\asymp n2^{-n}$ (see the proof of
Proposition~1.10 below). However, 0 is the only point
of the interval $[0,1]$, where the density $\frac{d\nu}{d\mu}$ is unbounded
(see Appendix~C).
\endremark
Let $R$ denote the rotation of the circle $\Bbb R/\Bbb Z$ by the
angle $\la^{-1}$.

\proclaim{Corollary 1.9} The Erd\"os measure $\mu$ is quasi-invariant
ergodic with respect to $T$ and $R$.
\endproclaim
\demo{Proof} The first claim follows directly from the equivalence
of the measures $\mu$ and $\nu$. To prove the second one, we note that
by Lemma~1.1,
$$
\mu(T^{-1}E)=\frac12\bigl(\mu E+\mu(E+\la^{-2}\mod 1)\bigr), \tag 1.7
$$
whence follows the required assertion.
\enddemo

We recall the following well-known claim which is a corollary of the
individual ergodic theorem. Namely,
two Borel measures invariant and ergodic with respect to one and the same
transformation of a metric space, either coincide or are mutually
singular.

Now we can present a new (dynamical) proof of the Erd\"os theorem on
the singularity of the measure $\mu$.

\proclaim{Proposition 1.10} {\rm (Erd\"os theorem, see \cite{Er})}
The Erd\"os measure
$\mu$ is singular with respect to Lebesgue measure $m$.
\endproclaim
\demo{Proof} We proved that $T\nu=\nu$, and above it
was noted that $Tm'=m'$ (see item~1.2),
hence by the corollary of the ergodic theorem,
either $\nu\perp m'$ or $\nu=m'$.
 Indeed, we can apply it, because $\nu$ is
ergodic by Proposition~1.6, and and the ergodicity of
$m'$ is a classical fact (see, e.g., \cite{Ge},
\cite{Re}). To show that $\nu\neq m'$,
we observe that $m'(0,\la^{-n})\asymp\la^{-n}$, while
$\nu(0,\la^{-n})=\wt\nu(\e_1=\e_2=\dots=\e_n=0)=O(n2^{-n})$,
because if for a sequence $x=\{x_n\}\in\wt\Sigma,\ \wt\n(x)\in
(\e_1=\e_2=\dots=\e_n=0)$, then either $x_i=0,\ 1\le i\le n$, or
$x_i=0,\ k+1\le i\le n$, and $x_k=1, x_{k-1}=1, x_{k-2}=0, x_{k-3}=1,
x_{k-4}=0$, etc. for some $k$. So, $\nu\perp m'$, hence, $\mu\perp m_1$,
as $m'\approx m_1,\ \mu\approx\nu$. \qed
\enddemo
\remark{Remark $1$} Note that the initial proof of Erd\"os followed
the traditions of those times and was based on the study of the Fourier
transform of $\mu$.
\endremark
\remark{Remark $2$} The present proof
fills a gap in the proof of this statement
in the previous joint paper by the authors \cite{VerSi}.
Another dynamical proof of Erd\"os theorem is given in Corollary~2.8
(see below).
\endremark

\remark{Remark $3$} The problem of computing the densities
$\frac{d\nu}{d\mu},\ \frac{d(R\mu)}{d\mu}$ and
$\frac{d(T\mu)}{d\mu}$ will be solved in Appendix~C.
Note that all these densities prove to be piecewise constant
and unbounded.
\endremark

We recall that $\wt m$ is the Markov measure on $\wt X$ with
maximal entropy (see item~1.4)

\proclaim{Proposition 1.11} The two-sided Erd\"os measure $\wt\nu$
is singular with respect to the Markov measure $\wt m$.
\endproclaim
\demo{Proof} We again apply the corollary of
the ergodic theorem to the transformation
$\wt\tau$ and the measures $\wt m$ and $\wt\nu$ on the two-sided
Markov compactum.
The distinction of the two measures is a consequence of the noncoincidence
of their one-sided restrictions (see Proposition~1.10), therefore,
$\wt m\perp\wt\nu$.
\enddemo
At the end of the item we prove a claim which we will need in the next
item. Recall that $\wt X$ has an additive structure (item~1.4).

\proclaim{Proposition 1.12} The Erd\"os measure $\wt\nu$ is invariant
under the transformation $i:\{\e_n\}\mapsto -\{\e_n\}$.
\endproclaim
\demo{Proof} Note that $i(\{\e_n\})=\wt\n(\{\e'_n\})$, where
$\{\e'_n\}\in\wt\Sigma$, and $\e'_n=1-\e_n$.
Therefore, for any Borel $E\subset \wt X$ and any sequence $\{x_n\}$
from the set
$\wt\n^{-1}(-E)$  there exists a unique sequence $\{x'_n\}\in\wt\n^{-1}(E)$
such that $x'_n=1-x_n$. Now the claim of the proposition follows from
the definition of $\wt\nu$ and the symmetricity of the measure $\wt p$
on $\wt\Sigma$.
\enddemo

\subhead 1.6. Erd\"os measure on the 2-torus \endsubhead
There exists an important smooth interpretation of the two-sided theory.
It is related to a general arithmetic approach to the coding of
the hyperbolic automorphisms of the torus.
Here we will give only some definitions and primary claims whose aim is
to describe a two-sided analog of the Erd\"os measure. Some necessary
bibliographic references will be given at the end of the item.

Consider the {\it Fibonacci automorphism} $\wt T$ of the 2-torus, i.e.
the automorphism with the matrix
$\(\smallmatrix 1&1\\1&0 \endsmallmatrix\)$. Later it will be clear that
this automorphism can be considered as a natural extension of
the endomorphism $Tx=\{\la x\}$ of the interval.

There exists a natural way to define a mapping semiconjugating the shift
$\wt\tau$ on the Markov compactum $\wt X$ and the Fibonacci automorphism,
namely the mapping which naturally generalizes
the $\beta$-expansions with $\beta=\la$
to the two-sided case. It is defined by the formula
$$
\wt l(\{\e_k\}_{-\infty}^\infty)=
\(\sum_{k=-\infty}^{\infty}\e_k\la^{-k}\mod 1,
\sum_{k=-\infty}^{\infty}\e_k\la^{-k-1}\mod 1\).\tag 1.8
$$
The convergence of the series involved follows from
$\la$ being a {\it PV number}.\footnote{I.e.
an algebraic integer greater than 1
whose Galois conjugates have the moduli less than 1, see, e.g.,
\cite{Cas}.} Indeed, as $-\la^{-1}$ is the Galois conjugate of $\la$,
we have $\|\la^n\|\le \la^{-n}$ for any $n\ge1$, where as usual,
$\|x\|:=\min\,(\{x\},1-\{x\})$.
Let us explain the background of formula~(1.8). Consider first
$x\ge0$ and its expansion $x=\sum_{k=-\infty}^\infty\e_k\la^{-k}$
which is the canonical expansion natural extended
to all nonnegative reals with $\e_k\equiv0$ for $k\le K(x)$.
So, we identify the set of sequences finite to the left with $\Bbb R_+$.
Consider now the inclusion
$\Bbb R_+\leftrightarrow \Cal R=
\{(\{x\},\{\la^{-1}x\})\mid x\ge0\}\subset \Bbb T^2$.
Since the set $\Cal R$ is the half-leaf of the unstable foliation
for the Fibonacci automorphism (corresponding to its eigenvalue
$\la$), we make sure that $(\wt l\wt\tau)(\{\e_n\})=\wt T\wt l(\{\e_n\})$,
where $\{\e_n\}$ is finite to the left.

As the
set $\Cal R$ is dense in the 2-torus, as well as the set of sequences finite
to the left is dense in $\wt X$, we can extend the relation above
to the whole compactum $\wt X$, i.e.
$$
\wt l\wt\tau=\wt T\wt l
$$
everywhere.

Besides, $\wt l$ is surjective and
from the proposition cited in item~1.4 and the fact that
Lebesgue measure $m_2$ is the only measure invariant under the
translations by a dense set of points of the 2-torus, it follows
that $m_2=\wt l(\wt m)$.

The important property of the mapping $\wt l$ is that
it is {\bf not} bijective a.e. Nevertheless, as is well known,
the automorphisms $(\Bbb T^2, m_2, \wt T)$ and $(\wt X, \wt m, \wt\tau)$
are metrically isomorphic, see, e.g., \cite{AdWe}. Below we will
introduce a conjugating mapping for these dynamical systems
which will be bijective a.e.

We recall that after a certain identification of pairs of sequences of
measure zero, the Markov compactum $\wt X$ becomes a group in addition
which we denoted by $\wt X'$ (see item~1.4). Note that the mapping $\wt l$
is well defined on $\wt X'$, i.e. it is constant on the equivalence classes.

\proclaim{Lemma 1.13} The mapping $\wt l:\wt X'\to\Bbb T^2$ is a
group homomorphism.
\endproclaim
\demo{Proof} It suffices to check that
$\wt l(e_{n-1})=\wt l(e_n)+\wt l(e_{n+1})$, where $e_k$
is a sequence having 1 at the $k$'th place and 0 at the other places.
This follows directly from formula~(1.8).
\enddemo

Now we are going to prove the following assertion.

\proclaim{Proposition 1.14} The mapping $\wt l$ is 5-to-1 a.e.
\endproclaim
We are going to give two different proofs of this proposition.
\demo{First (geometric) proof}
Consider an arbitrary sequence $\e=\{\e_n\}_{-\infty}^\infty\in\wt X$.
We split it into two pieces $\{\e_n\}_{-\infty}^0$ and $\{\e_n\}_1^\infty$
and define $x_1(\e):=\sum_{k=1}^\infty \e_k\la^{-k},\ x_2=\sum_{k=0}^\infty
\e_{-k}(-\la)^{-k}$. It is a direct inspection that $x_1\in[0,1],\
x_2\in[-1,\la]$. Using the relation $\{\la^n\}=\{(-1)^{n+1}\la^{-n}\},\
n\ge0$, we make sure that $\sum_{-\infty}^\infty\e_n\la^{-n}=x_1-x_2\mod1$
and similarly,
$\sum_{-\infty}^\infty\e_n\la^{-n-1}=\la^{-1}x_1+\la x_2\mod1$.
Thus, we have the sequence of mappings
$$
\wt X @>\phi>> \Bbb R^2 @>b>> \Bbb R^2,
$$
where $\phi(\e)=(x_1,x_2)$, and $b(x_1,x_2)=(x_1-x_2,\,\la^{-1}x_1+\la x_2)$.
Now
$$
\wt l(\e)=(b\phi)(\e)\mod\Bbb Z^2.
$$
Note that since $(\e_0,\e_1)\neq(1,1)$,
the $\phi$-image of $\wt X$ is in fact the difference of
the rectangles
$\Pi=([0,1]\times[-1,\la])\setminus([\la^{-1},1]\times[\la^{-1},\la])$.
As the area of $\Pi$ is $\la^2-\la^{-2}=\sqrt5$, and the linear
transformation $b=\(\smallmatrix 1 & -1 \\ \la^{-1} & \la\endsmallmatrix\)$
from $\Bbb R^2$ to $\Bbb R^2$ has determinant $\sqrt5$, the
image $(b\phi)(\wt X)$ has area 5. Since this mentioned image is
also the difference of some rectangles whose vertices are
$(1,-\la),\ (2,-1),\ (\la^{-2},\la),\ (0,\la^{-1}\sqrt5),\ (-1,3),\
(-\la,\la^2)$, one can immediately check that this set is
a 5-tuple fundamental domain with respect to the lattice $\Bbb Z^2$.\qed
\enddemo
\remark{Remark} On the other hand, the final statement of the
proof follows from Lemma~1.13, as from its claim it follows that
the $\wt l$-preimage of a.e. point of the 2-torus has one and the same
capacity.
\endremark
\demo{Second (algebraic) proof}
From Lemma~1.13 it follows that
one needs only to describe the kernel of the homomorphism $\wt l$.
From the general considerations we conclude that the mapping
$\wt l$ is bounded-to-one, whence $\#\Ker\wt l<+\infty$.
Furthermore, this kernel is obviously invariant under the shift $\wt\tau$,
hence it consists of purely periodic sequences only.

Let $\e=(\e_1,\dots,\e_r)^\infty$ be such a sequence. From formula~(1.8)  
it follows that
$$
\|\xi\la^n\|\to0,\quad n\to +\infty,
$$
where $\xi=\xi(\e)=(\sum_1^r\e_k\la^{-k})/(\la^r-1)\in\Bbb Q(\la)$.
Let $\Cal G:=\{\xi:\|\xi\la^n\|\to0,\ n\to\infty\}$. Obviously,
$\Cal G$ is a group in addition.
The following lemma answers the question on the structure of the
group $\Cal G$ which is important itself and will be used
below. Let $\bold Z[\la]$ denote the additive group of the ring of
all Laurent polynomials in powers of $\la$,
i.e. $\{m+n\la\mid m,n\in\Bbb Z\}$.

\proclaim{Lemma 1.15} The group $\Cal G$ is isomorphic to $\Bbb Z^2$.
Its elements are described as follows:
$$
\Cal G=\left\{\frac{m+n\la}5\, :\,m,n\in\Bbb Z,\ 2n-m\equiv0\pmod5\right\}.
$$
The factor group $\Cal H=\Cal G/\bold Z[\la]$ is the cyclic group
$\{a\frac{\la+2}5\mid a\in\Bbb Z/5\Bbb Z\}$.
\endproclaim
\demo{Proof}
By the theorem of Pisot-Vijayaraghavan
on the structure of the group $\Cal G$
(see \cite{Cas}), $\xi$ is necessarily algebraic and belongs
to the field $\Bbb Q(\la)$, and
the necessary and sufficient condition for $\xi$ to belong to $\Cal G$
is $\Tr(\xi)\in\Bbb Z, \Tr(\la\xi)\in\Bbb Z$, where $\Tr$ denotes the trace
of an algebraic number. Solving these equations for $\xi=(m+n\la)/q$
with $m,n\in\Bbb Z,\ q\in\Bbb N$ and at least one of the numbers
$m,\ n$ being coprime with $q$, we come to the system of congruences
$$
\cases
2m+n\equiv0 & \pmod q \\
m+3n\equiv0 & \pmod q,
\endcases
$$
whence $5m\equiv0\pmod q,\ 5n\equiv0\pmod q$, and thus, $q=1$ or 5.
Besides, the system above with $q=5$ is equivalent to
one congruence $2n-m\equiv0\pmod5$. The second claim of the lemma
is a direct computation.
\enddemo

Return to the second proof of the proposition.
From the formula for $\xi=\xi(\e)$ above it follows that $\xi$ belongs to
$\bold Z[\la]$ if and only if its period $r\le2$, i.e. when the sequence
$\e$
is equivalent to $0^\infty$ in $\wt X'$ (see the definition of
$\wt X'$ in item~1.4). Besides, if for two periodic sequences $\e$ and $\e'$,
$\xi(\e)$ and $\xi(\e')$ belong to one and the same element of the
factor group $\Cal H$, we see that $\wt l(\e)=\wt l(\e')$.

Thus, $\#\Ker\wt l\le5$. Let
the sequence $\e^{(0)}=(0.100)^\infty$; consider the set
$\Cal H'=\{0,\ \e^{(j)},\ 0\le j\le 3\}$, where
$\e^{(j)}=\wt\tau^j(\e^{(0)})$.

It is verified directly that $\Cal H'$ is a subgroup
of $\wt X'$. We claim that $\Cal H'=\Ker\wt l$. By the above, it suffices to
prove
the inclusion $\Cal H'\subset\Ker\wt l$. From formula~(1.8)  ,
$\wt l(\e^{(j)})=(a,a),\ 0\le j\le 3$, where
$$
a=\lim_{n\to\infty}\left\|\frac{\la^n}{\la^4-1}\right\|=
\lim_{n\to\infty}\left\|\frac{\la^n}{\la^2+1}\right\|=
\lim_{n\to\infty}\left\|\frac{\la^n}{\sqrt5}\right\|=0,\tag 1.9
$$
as $\la^n=F_{n-1}\la+F_{n-2},\ n\ge 3$, and we have
$\la^n\sim F_{n-1}\sqrt5$, whence $\|\frac1{\sqrt5}\la^n\|\to0$
as $n\to\infty$.
\enddemo
\remark{Remark} There exists a simple explanation of the origin
of the sequences $\e^{(k)}$. Note first
that if we identify a positive integer $n$ with the sequence
equal to $ne_0$, we obtain a natural inclusion $\Bbb N\subset\wt X$.
For example, $2=\la+\la^{-2},\ 3=\la^2+\la^{-2}$, etc.
It is a direct computation that
the Fibonacci numbers have the following representations in the
compactum $\wt X$:
$$
\align
F_k&=\la^{k-1}+\la^{k-5}+\cdots+\la^{-k+3}+\la^{-k}, \quad k \text{ even},\\
F_k&=\la^{k-1}+\la^{k-5}+\cdots+\la^{-k+5}+\la^{-k+1}, \quad k \text{ odd}.
\endalign
$$
Thus, $\e^{(j)}=\lim_k(F_{4k+j}e_0),\,0\le j\le3$ in the compactum $\wt X$.
Since $\lim_k\|\la F_k\|=0$, we have
$\wt l(\e^{(j)})=(0,0)$.
\endremark

Now we are going to present
a simple modification of $\wt l$ which proves to be a bijection.
\definition{Definition} Let the mapping $\wt L:\wt X\to \Bbb T^2$ be defined
by the formula
$$
\wt L(\{\e_k\}_{-\infty}^\infty)=
\(\sum_{k=-\infty}^{\infty}\e_k\frac{\la^{-k}}{\sqrt5}\mod 1,
\sum_{k=-\infty}^{\infty}\e_k\frac{\la^{-k-1}}{\sqrt5}\mod 1\). \tag 1.10
$$
\enddefinition
By formula~(1.9) , a
series $\sum_{k\in\Bbb Z}\e_k\frac{\la^{-k}}{\sqrt5}$ converges
modulo 1 for any $\{\e_k\}\in\wt X$.
\remark{Remark} In fact, we can treat formula~(1.10) as formula~(1.8) with
the set of digits $\{0,\ 1/\sqrt5\}$ instead of $\{0,1\}$ (see
the remark about the references at the end of the item).
\endremark

Furthermore, $\wt L$ semiconjugates the shift and the Fibonacci automorphism,
and by the same purposes as for $\wt l$,
we have $m_2=\wt L(\wt m)$.

\proclaim{Proposition 1.16} The mapping $\wt L:\wt X\to \Bbb T^2$ is
bijective a.e.
\endproclaim
\demo{Proof}
Let the projection
$P:\wt X\to\wt X$ be defined by the formula $P(e_n)=e_{n-1}+e_{n+1}$
and extended to the whole compactum $\wt X$ by linearity (here
$e_n$ is, as above, the sequence having the only 1 at the $n$'th place).
Then by the relation $\la^k=\la^{k-1}/\sqrt5+\la^{k+1}/\sqrt5$
and formulas~(1.8)  and (1.10) , we have $\wt l=\wt L P$. Now, considering
$A=\wt L P(\wt L)^{-1}:\Bbb T^2\to\Bbb T^2$, we see that
$A=\wt T+\wt T^{-1}=\(\smallmatrix 1&2\\2&-1\endsmallmatrix\)$.
Since $|\det(A)|=5$, $\wt l$ is 5-to-1 and $\wt l=\wt L A$,
we complete the proof.
\enddemo
So, we proved the following theorem.

\proclaim{Theorem 1.17} The mapping $\wt L$ is a metric isomorphism of
the two-sided Markov shift $(\wt X,\wt m,\wt\tau)$ and the Fibonacci
automorphism $(\Bbb T^2,m_2,\wt T)$.
\endproclaim
\remark{Remark {\rm1}} Furthermore, $\wt L$ is a group
isomorpism of the groups $\wt X'$ and $\Bbb T^2$.
\endremark
\remark{Remark {\rm2}} Actually, the mappings $\wt l$ and $\wt L$ lead
to different interpretations of the torus as the image of $\wt X$. Namely,
$\Bbb T^2\cong(\Bbb Q(\la)\otimes\Bbb R)/(\Bbb Z+\la\Bbb Z)$ for $\wt l$,
while for $\wt L$ we have
$\Bbb T^2\cong(\Bbb Q(\la)\otimes\Bbb R)
\left/\frac{\Bbb Z+\la\Bbb Z}{\sqrt5}\right.$.
\endremark

Now we are ready to define the two-dimensional Erd\"os measure.

\definition{Definition}
Let the measure $\wt\mu$ on the 2-torus be defined as
$\wt\mu=\wt L(\wt\nu)$. We call it the {\it two-dimensional Erd\"os}
measure.
\enddefinition

Let us show that in a sense the two-dimensional Erd\"os measure is defined
canonically. Note that any mapping $\Psi$
from $\wt X$ onto $\Bbb T^2$ such that
\roster
\item $\Psi\wt\tau=\wt T\Psi$,
\item $\Psi(\e+\e')=\Psi(\e)+\Psi(\e')$ for a.e. $\e,\ \e'\in\wt X$,
\item $\Psi$ is one-to-one a.e.
\endroster
is of the form
$$
\Psi(\{\e_k\}_{-\infty}^\infty)=\wt l_\xi(\{\e_k\}_{-\infty}^\infty)=
\(\sum_{k=-\infty}^\infty\e_k\xi\la^{-k}\mod 1,
\sum_{k=-\infty}^\infty\e_k\xi\la^{-k-1}\mod 1\),\tag 1.11
$$
where $\xi$ is some real number.
\remark{Remark} Dealing with the leaf of the {\bf stable} foliation for
$\wt T$ (see the explanation after formula~(1.8) above), we come to
a similar family of mappings:
$$
\wt\kappa_\xi(\{\e_k\}_{-\infty}^\infty)=
\(\sum_{k=-\infty}^\infty\e_k\xi(-\la)^k\mod 1,
\sum_{k=-\infty}^\infty\e_k\xi(-\la)^{k+1}\mod 1\).
$$
However, it does not yields anything new, since below we will
show that $\xi$ should be a quadratic irrational,
and it is easy to see that for such a $\xi$,
$\wt\kappa_\xi=\wt l_{-\ov\xi}$, where $\ov\xi$ denotes
the algebraic conjugate of $\xi$.
\endremark
We formulate
the claim which answers the question on the possible
values for $\xi$ and the bijectivity of
$\wt l_\xi$.

\proclaim{Proposition 1.18}
\roster
\item For the series in formula~(1.11) to converge,
$\xi$ must belong to the group $\Cal G$.
\item The mapping $\wt l_\xi$
is one-to-one
a.e. if and only if $\xi=\pm\frac{\la^k}{\sqrt5}$ for $k\in\Bbb Z$.
\endroster
\endproclaim
\demo{Proof} The first claim is a consequence of the condition
$\|\xi\la^n\|\to0,\ n\to \infty$,
which defines the group $\Cal G$. Let us prove the second one.
Using the same geometric arguments and the same notation as in the
proof of Proposition~1.14 and also the fact that
$\ov{\xi\la^n}=\ov\xi(-\la)^{-n}$, we
obtain the following sequence of mappings
$$
\wt X @>\phi>> \Pi @>b_\xi>> \Bbb R^2,
$$
where
$$
b_\xi=\pmatrix \xi & -\ov\xi \\ \la^{-1}\xi & \la\ov\xi \endpmatrix,
$$
and $(\wt l_\xi)(\e)=(b_\xi\phi)(\e)\mod\Bbb Z^2$.
The area of the set $(b_\xi\phi)(\wt X)$ is $5|\xi\ov\xi|=5|N(\xi)|$
(the algebraic norm of $\xi$).
So, the condition of the bijectivity of $\wt l$ a.e. is
$$
N(\xi)=\pm\frac15,
$$
where $\xi\in\Cal G$. We are going to use Lemma~1.15; let
$\xi=\frac{m+n\la}5$ with the extra condition
$$
2n-m\equiv0\pmod5.\tag *
$$
Then
$N(\xi)=\frac1{25}(m^2+mn-n^2)$, and we come to the Diophantine equation
$$
m^2+mn-n^2=\pm5
$$
together with (*). Putting $u=n,\ v=\frac{2m+n}5\in\Bbb Z$,
we get a classical Pell's equation
$$
u^2-5v^2=\pm4.
$$
Using the usual method (see, e.g., the monograph
\cite{Lev, vol.~1, Th.~8--7}), we obtain its general solution in the form
$(u,v)=(u_k,v_k)$, where
$$
u_k+\sqrt5 v_k=\pm2\la^k,\quad k\in\Bbb Z.
$$
Returning to the variables $m=\frac{5v-u}2$ and $n=u$, we have
$$
\frac{2n-m + (2m+n)\la}5=\pm\la^k,\quad k\in\Bbb Z,
$$
whence
$$
\xi=\pm\frac{\la^k}{\sqrt5},\quad k\in\Bbb Z.\qed
$$
\enddemo
\remark{Remark {\rm1}} Note that the condition~(*) proved to be satisfied
automatically.
\endremark
\remark{Remark {\rm2}} Thus, any value of a parameter $\xi$ yielding
the bijectivity of $\wt l_\xi$ is of the form $\frac1{\sqrt5}$ times
an arbitrary unit of the field $\Bbb Q(\la)$.
\endremark

\proclaim{Corollary 1.19} For any bijection $\wt l_\xi$
from $\wt X$ onto $\Bbb T^2$ having the form~(1.11),
the image of the Erd\"os measure $\wt\nu$ under it will coincide
with $\wt\mu$.
\endproclaim
\demo{Proof} It suffices to use the form of any mapping of such a kind
deduced in
the previous proposition and the invariance of $\wt\nu$ under the
shift $\wt\tau$ (Proposition~1.3) and the operation $\e\mapsto -\e$
(Proposition~1.12).
\enddemo

So, the two-dimensional Erd\"os measure is defined canonically,
and we show that the result analogous to the Erd\"os theorem
holds for the two-dimensional Erd\"os measure.

\proclaim{Theorem 1.20}
The two-dimensional Erd\"os measure $\wt\mu$
on $\Bbb T^2$
is singular with respect to Lebesgue measure.
\endproclaim
\demo{Proof} It follows from the Proposition~1.11 that $\wt m$
is singular with respect to $\wt\nu$,
and by the bijectivity of $\wt L$ a.e. with respect to both measures, we
deduce the mutual singularity of their images.
\enddemo

\remark{Remark} Let us formulate a number of open questions about the
properties of the two-dimensional Erd\"os measure.
\roster
\item Is the measure $\wt\mu$ Gibbs with respect to the Fibonacci
automorphism for some natural potential?
\item Is $\wt\mu$ invariant under the endomorphism $A=
\(\smallmatrix 1&2\\2&-1\endsmallmatrix\)$ (if it was so, the mapping
$\wt l(\wt\nu)$ would coincide with $\wt\mu$)?
\item Both questions can be reformulated in terms of the compactum
$\wt X$ and remain valid for it.
\endroster
\endremark

We conclude this item by mentioning some necessary references.
The first precise symbolic coding of the hyperbolic automorphisms
of the 2-torus had been proposed in \cite{AdWe} and was developing
then by a number of authors (see the references in \cite{Ver6},
\cite{KenVer}). In \cite{Ber} in connection with the arithmetic of PV
numbers were considered two-sided expansions and the corresponding
mapping semicongugating the two-sided (in general, {\it sofic})
shift and the endomorphism of the torus with a companion matrix.
Note that for the case of the golden ratio it coincides with the
mapping $\wt l$.

In the works \cite{Ver3}, \cite{Ver5}, \cite{Ver6}, \cite{KenVer}
an arithmetic approach to the coding of hyperboilc automorphisms
of the torus has been developing. In one of the versions of such an
approach which generalizes
adic transformation to the two-sided case, it leads to
the two-sided $\beta$-expansions, and another one being applicable
to a more general algebraic numbers (not only PV) leads to
a scheme of $\wt L$; it uses the digits from the field
of an irrationality being not always integers
(see \cite{Ver6}, \cite{KenVer}).
Later this approach was developed and detailed in the dissertation
\cite{Leb}.

The group $\Cal H'$ and its action on $\wt X$ (see above) were considered
in the recent work \cite{FrSa} in connection with the study of certain
finite automata.

More detailed analysis of bijective arithmetic codings for hyperbolic
automorphisms of the 2-torus was recently given in \cite{SiVer}.

\subhead 1.7. Numerical properties of the measure $\nu$ \endsubhead
We are going first to give now an explicit formula for $\nu$.

\proclaim{Proposition 1.21} The following relation holds:
$$
\nu E=
\cases \frac23\mu E+\frac13\mu(E+\la^{-2})+\frac16\mu(E+\la^{-1}),
       & E\subset [0, \la^{-2}) \\
             \frac23\mu E+\frac13\mu(E+\la^{-2}),
       & E\subset [\la^{-2}, \la^{-1}) \\
             \frac12\mu E+\frac13\mu(E-\la^{-1}),
       & E\subset [\la^{-1}, 1].
\endcases\tag 1.12
$$
\endproclaim
\demo{Proof} Let the measure $\nu'$ be defined by formula~(1.12). We need
to show that $\nu'=\nu$. Note first that the
$T$-invariance of $\nu'$ is checked directly using Lemma~1.1.
Let, say, $E\subset (0,\la^{-2})$. Then $T^{-1}E=\la^{-1}E\cup
(\la^{-1}E+\la^{-1})$. Hence
$$
\align
\nu'(T^{-1}E)&= \nu'(\la^{-1}E)+\nu'(\la^{-1}E+\la^{-1})\\
                  &= \frac23\mu(\la^{-1}E)+\frac13\mu(\la^{-1}E+\la^{-2})
                     +\frac16\mu(\la^{-1}E+\la^{-1})+
                     \frac12\mu(\la^{-1}E+\la^{-1})\\
                  &+ \frac13\mu(\la^{-1}E)\\
                  &= \frac12\mu E+\frac13\mu(E+\la^{-2})+
                     \frac16\mu(E+\la^{-1})+\frac16\mu E
                     \quad\text{(by Lemma 1.1)}\\
                  &= \frac23\mu E+\frac13\mu(E+\la^{-2})+
                     \frac16\mu(E+\la^{-1}) \\
                  &= \nu'(E).
\endalign
$$
The cases $E\subset (\la^{-2},\la^{-1})$ and $E\subset (\la^{-1},1)$
are studied in the same way.

To prove now that $\nu'=\nu$, we observe that since $\mu$ is quasi-invariant
with respect to $T$ and the rotations by $\la^{-1}$ and $\la^{-2}$,
the measure $\nu'$ is also ergodic with respect to $T$ (see (1.7) and
(1.12)). Since $\mu\approx\nu$ and $\mu\prec\nu'$, we have $\nu\prec\nu'$,
and by the corollary of the ergodic theorem, $\nu'=\nu$.\qed
\enddemo

\proclaim{Corollary 1.22}
$\nu(0,\la^{-2})=\frac49,\ \nu(\la^{-2},\la^{-1})=\nu(\la^{-1},1)=\frac5{18}$.
\endproclaim

\proclaim{Proposition 1.23}
The following relation holds:
$$
\nu=\lim_n T^n\mu.
$$
\endproclaim
\demo{Proof}
The sketch of
the proof is as follows. Let, for example, $E\subset(0,\la^{-2})$.
Considering successively the sets $T^{-n}E,\ n\ge1$ and
using Lemma~1.1, we deduce similarly to relation~(1.7) that
$\mu(T^{-2}E)=\frac34\mu E+\frac14\mu(E+\la^{-2})+\frac14\mu(E+\la^{-1}),\
\mu(T^{-3}E)=\frac58\mu E+\frac38\mu(E+\la^{-2})+\frac18\mu(E+\la^{-1})$,
etc., whence by induction, $\mu(T^{-n}E)=\(\frac12+\frac14-\frac18+\dots+
\frac{(-1)^n}{2^n}\)\mu E \allowmathbreak
+ \(\frac12-\frac14+\dots+\frac{(-1)^{n+1}}{2^n}\)
\mu(E+\la^{-2}) +
\(\frac14-\frac18+\dots+\frac{(-1)^n}{2^n}\)\mu(E+\la^{-1})=
\frac23\mu E + \frac13\mu(E+\la^{-2})+\frac16\mu(E+\la^{-1})+O(2^{-n})=
\nu E+O(2^{-n})$ by formula~(1.12). The cases $E\subset(\la^{-2},\la^{-1})$
and $E\subset(\la^{-1},1)$ are considered in the same way.
\enddemo

\remark{Remark} It is appropriate,
following the well-known framework of the baker's transformation
which serves as a model for the full two-sided shift on $\wt\Sigma$,
to represent the two-sided shift on
$\wt X$ as the
{\it Fibonacci-baker's} transformation.

Namely, we split a sequence $\{\e_k\}\in \wt X$ into the two one-sided
sequences, i.e. into $(\e_1\e_2\dots)\in X$ and $(\e_0\e_{-1}\dots)\in X$
with regard to the fact that $\e_0\e_{-1}=0$. This last condition leads
to the space $Y=([0,1]\times[0,1])\setminus([\la^{-1},1]\times[\la^{-1},1])$
similar to the set $\Pi$ described in the proof of Proposition~1.14.

\bigskip
\epsfysize=4cm
\centerline{\epsfbox{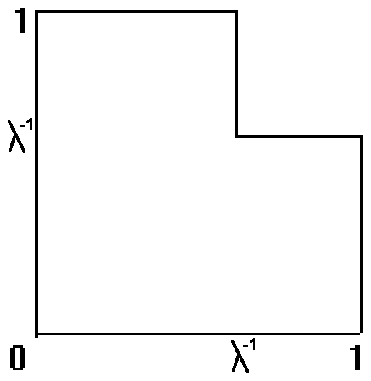}}
\bigskip\bigskip
\nopagebreak

\centerline{Fig. 1. The natural domain for the Fibonacci-baker's
transformation}
\medskip

Thus, the shift $\wt\tau$ on the two-sided Markov compactum $\wt X$
is isomorphic to the transformation $F$ on the space $Y$ with
$$
F(x,y)=\cases \bigl(\la x,\la^{-1}y\bigr),&          x\in[0,\la^{-1}] \\
              \bigl(\la x-1,\la^{-1}y+\la^{-1}\bigr),& x\in(\la^{-1},1].
       \endcases
$$
We call $F$ the {\it Fibonacci-baker's} transformation on the set $Y$
(see Fig.~1). For more general models this transformation was considered
in \cite{DaKrSo}.
\endremark

\head 2. Symbolic dynamics of expansions \endhead

In this section we will study in detail the combinatorics of all
possible representations of a real $x$ of the form~(1.1) with
$\e_k\in\{0,1\}$ for all $k$.

\subhead 2.1. Blocks \endsubhead Let us give an important
technical definition.
\definition{Definition} A finite 0-1 sequence
without pairs of adjacent 1's starting
from 1 and ending by an even number of zeroes will be called
a {\it block} if it does not contain any piece ``$1(00)^l1$''
with $l\ge 1$.
\enddefinition

Let us make some remarks. Note first that each block has odd length;
the simplest example of a block is ``100''. Next, there are exactly
$2^{n-1}$ blocks of length $2n+1$. This assertion follows from
the fact that a block $B$ can be represented in the form
$1(00)^{a_1}(01)^{a_2}(00)^{a_3}\dots (01)^{a_{t-1}}(00)^{a_t}$ for $t$ odd
or $1(01)^{a_1}(00)^{a_2}\dots (01)^{a_{t-1}}(00)^{a_t}$ for
$t$ even.
Thus, any block $B$ is naturally parametrized by means of a finite
sequence of positive integers $a_1,\dots, a_t$,
and we will write $B=B(a_1,\dots,a_t)$.

Let $I_0:=[\la^{-1}, 1)$, i.e. the interval corresponding to the cylinder
$(\e_1=1)\subset X$.
\definition{Definition}
Let $x$ lie in the interval $I_0$, and let
the canonical expansion of $x$ have infinitely many pieces
``$1(00)^l1$'' with $l\ge1$. We call such a point $x$ {\it regular}.
\enddefinition
Almost every point $x$ in $I_0$ with respect
to Lebesgue measure is regular.
Now we split the canonical expansion of a regular $x$ into blocks as follows.
Since $x\in I_0$, its canonical expansion
starts with 1. It is just the beginning of the first block
$B_1=B_1(x)$. The first block ends, when an even number of zeroes
followed by 1 appears for the first time. This 1 begins
the second block $B_2=B_2(x)$ of the canonical expansion of $x$, etc.
We define thus a one-to-one mapping $\Psi$ acting from the set of all
regular points of $(\la^{-1},1)$ into the space of block
sequences.
\definition{Definition} The sequence $(B_1(x),B_2(x),\dots)=\Psi(x)$
will be called the {\it block expansion} of a regular $x$.
\enddefinition

\subhead 2.2. The cardinality of a 0-1 sequence and its properties
\endsubhead
We are going to define an equivalence relation on the set of all
finite 0-1 sequnces.
\definition{Definition} Two 0-1 sequences
(finite or not) $(x_1x_2\dots)$ and
$(x'_1x'_2\dots)$ are called {\it equivalent} if $\sum_k x_k\la^{-k}=
\sum_k x'_k\la^{-k}$ (or, equivalently, if their normalizations coincide
--- see Section~1).
Let for a finite 0-1 sequence $x$,
$\E(x)$ denote the set of all 0-1 sequences equivalent
to $x$; this set is always finite.
Let $f(x):=\#\E(x)$. We call $f$ the {\it cardinality}
of a finite sequence (or the cardinality of an equivalence class).
\enddefinition
Note that this function (of positive integers) was considered in \cite{Car},
\cite{AlZa} and recently in \cite{DuSiTh} and \cite{Pu}.

The assertions below answer the question about the cardinality of a block
and explain the purpose of the introduction of blocks as natural structural
units in this theory.
\proclaim{Lemma 2.1} Let a finite word $x$ be a block, i.e.
$x=1(00)^{a_1}(01)^{a_2}(00)^{a_3}\dots (01)^{a_{t-1}}(00)^{a_t}$
or $1(01)^{a_1}(00)^{a_2}\dots (01)^{a_{t-1}}(00)^{a_t}$.
Let $p/q=[a_1, \dots, a_t]$ be a finite continued
fraction. Then
$$
f(x)=p+q.
$$
\endproclaim
\demo{Proof} Note first that $f(100)=2=p+q$. The desired relations
for the blocks 10000 and 10100 follow by direct inspection. Next, let
$\ka_k=\ka_k(a_1,\dots, a_k)=f(x)$. We need to show that
similarly to the numerators and denominators of the convergents,
$\ka_t=a_t\ka_{t-1}+\ka_{t-2}$, whence the required assertion will follow.

Let, say, $t$ be odd and
$B=1(00)^{a_1}\dots (00)^{a_{t-2}}(01)^{a_{t-1}}(00)^{a_t}$.
We first consider all
0-1 sequences equivalent to $B$ and ending by
$(00)^{a_t}$. Their number is obviously
$\ka_{t-2}$, as they in fact should end by
$(01)^{a_{t-1}}(00)^{a_t}$.
Now we consider 0-1 sequences equivalent to $B$ but not ending by
$(00)^{a_t}$ and will show that their number is $a_t\ka_{t-1}$. Namely,
let
$$
\align
&B^{(1)}=1(00)^{a_1}\dots (00)^{a_{t-2}}(01)^{a_{t-1}-1}0011(00)^{a_t-1},\\
&B^{(2)}=1(00)^{a_1}\dots (00)^{a_{t-2}}(01)^{a_{t-1}-1}0(01)^21(00)^{a_t-2},
\\
&\dots \\
&B^{(a_t)}=1(00)^{a_1}\dots (00)^{a_{t-2}}(01)^{a_{t-1}-1}0(01)^{a_t}1.
\endalign
$$
It is clear that any $B^{(j)}$ is equivalent to $B$; now we observe
that the number of 0-1 words equivalent to $B$ and ending by
$(00)^{a_t-j},\ 1\le j\le a_t$, is exactly $\ka_{t-1}$, as the
replaceable part of $B^{(j)}$ with the fixed end $(00)^{a_t-j}$ is in fact
$1(00)^{a_1}\dots (00)^{a_{t-2}}(01)^{a_{t-1}-1}00$,
hence $\ka_t=a_t\ka_{t-1}+\ka_{t-2}$, as $[a_1,\dots, a_{t-1}-1, 1]=
[a_1,\dots, a_{t-1}]$. The case of even $t$ is studied in the same
way. \qed
\enddemo
\remark{Remark} To any rational $r\in (0,1)$  exactly
two blocks correspond, namely
with $r=[a_1,\dots,a_t]=[a_1,\dots,a_{t-1}, a_t-1,1]$, and the unique
block ``100'' corresponds to $r=1$.
\endremark
If $E_1$ and $E_2$ are two sets of sequences,
then henceforward $E=E_1E_2$ is the
concatenation of these two sets, i.e. any sequence in $E$ begins
with a word from $E_1$ and ends with a word from $E_2$.

\proclaim{Lemma 2.2} For any blocks $B_1,\dots, B_k$,
\roster
\item $\E(B_1\dots B_k)=\E(B_1)\dots\E(B_k)$.
\item The cardinality is blockwise multiplicative, i.e.
$$
f(B_1\dots B_k)=\prod_1^k f(B_i).
$$
\endroster
\endproclaim
\demo{Proof}  It suffices to prove item~(1).
Let us restrict ourselves to the case $k=2$ (the general
case is studied in the same way). We will see that there is no sequence
in $\E(B_1B_2)$ containing a triple 011 or 100 which crosses the ``border''
between the first $|B_1|$ digits and the last $|B_2|$. In other words,
we will show that
any sequence equivalent to $B_1B_2$ can be constructed as the concatenation
of a sequence equivalent to $B_1$ and a sequence equivalent to $B_2$.

Let ``$|$'' below denote the border in question. First, any sequence from
$\E(B_2)$ must begin either with 10 or with 01, hence, the situation
$(1|00)$ or $(0|11)$ is impossible. Next, a sequence in $\E(B_1)$ ends
in either 00 or 11 (see the proof of the previous lemma), neither
leading to $(10|0)$ or $(01|1)$. \qed
\enddemo
This simple result shows that the space of all equivalent infinite
0-1 sequences for a given regular $x$ splits into the direct infinite product
of spaces, the $k$'th space consisting of all finite sequences equivalent
to the block $B_k(x)$. So, we see that the notion of block, initially
arising in terms of the canonical expansion, can be naturally
extended to all representations. Below we will explain the geometric
sense of a block in terms of the Fibonacci graph (see p.~3.1).
\remark{Remark}
Note that this block partition appeared for the first time in \cite{Pu}
in other terms and for algebraic and combinatorial purposes.
Namely, let the partial ordering on a
space $\E(x)$ for some finite word $x$ be defined as follows.
We set $x\prec x'$ if there exists $k\ge2$ such that $x_{k-1}=0, x_k=1,
x_{k+1}=1, x'_{k-1}=1, x'_k=0, x'_{k+1}=0$, and $x_j=x'_j,\ |k-j|\ge2$.
 Next, one extends this ordering
by transitivity. It was shown in \cite{Pu}
that any equivalence class has the structure
of a distributive lattice in the sense of this order.
\endremark

\subhead 2.3. Goldenshift \endsubhead We are going to give one of the
central definitions of the present paper.

\definition{Definition} The transformation $\S$ acting from the set
of regular points of the interval $(\la^{-1},1)$ into itself by the formula
$$
\S x=\tau^{n(x)}x, \quad x \,\,\,\text{is regular,}
$$
where $n(x)$ is the length of the first block of the block
expansion of $x$, is called the {\it goldenshift}.
\enddefinition
\remark{Remark {\rm1}} The transformation $\S$ is piecewise linear. More
precisely,  if $(\la^{-1},1)=\bigcup_r \De_r$ is the partition
of $(\la^{-1},1)$ mod 0 into intervals corresponding to that of
$\B$ into the states of the first block, then $\S$ is linear inside
$\De_r=:[\a_r,\beta_r)$, and $\S(\a_r)=\la^{-1},\,\,\S(\beta_r)=1$.
\endremark
\remark{Remark {\rm2}} The transformation
$\S$ is a {\it generalized
power} of $\tau$ in the sense of Dye (see, e.g., \cite{Bel}). In other words,
the goldenshift is a random power of the
ordinary shift, as the number of shifted coordinates
depends on the length of the first block.
Note that the goldenshift is not an induced
endomorphism for $\tau$ but for the two-sided case it is (see
Proposition~2.3 and Theorem~2.12 below).
\endremark
\remark{Remark {\rm3}} Let $\X$ denote the space of block
sequences. The goldenshift may be treated as
the one-sided shift in the space $\X$, i.e.
$\S(B_1B_2B_3\dots)=(B_2B_3\dots)$.
\endremark
Now we are going to define the two-sided goldenshift.
Let, as above,  $\wt\tau$ denote the two-sided shift on $\wt X$.
In order to define the two-sided goldenshift  on the Markov compactum $\wt X$,
we give the following definition.
\definition{Definition} The $\wt\tau$-invariant set
$\wt X^{\text{reg}}\subset \wt X$
is defined as the set consisting of all sequences containing
pieces ``$10^{2l}1$'' with $l\ge1$ infinitely many times both to the left
and to the right with respect to the first coordinate.
\enddefinition
Obviously,
$\wt\nu\bigl(\wt X^{\text{reg}}\bigr)=1$, as by Proposition~1.4, the measure
$\wt\nu$ of any cylinder in $\wt X$ is positive, and by Proposition~1.5,
the automomorphism $(\wt X,\wt\nu,\wt\tau)$ is ergodic, so,
it suffices to apply the ergodic theorem.
Let $\wt X_0:=\bigcup_{k=1}^\infty (x_{-2k}=1,\,x_{-2k+1}=\dots=x_0=0,\,
x_1=1)$, and let $\wt X_0^{\text{reg}}=\wt X^{\text{reg}}\cap
\wt X_0$.
\definition{Definition}
The two-sided goldenshift
$\wt\S:\wt X_0^{\text{reg}}\to\wt X_0^{\text{reg}}$
on the Markov compactum $\wt X$ is, by definition,
the shift by the length of $B_1$.
\enddefinition
\remark{Remark} Considering the space of two-sided
block sequences $\wt\X=
\prod_{-\infty}^\infty \B$ and implying that
$B_1$ begins with the first coordinate
of $\wt X$ (i.e. with $\e_1=1$), we  see that the two-sided goldenshift
$\wt\S$ is the shift in the space $\wt\X$.
\endremark

\proclaim{Proposition 2.3} The two-sided shift $\wt\tau$ on the set
$\wt X^{\text{reg}}$
is a special
automorphism over the goldenshift $\wt\S$ on the set $\wt X_0^{\text{reg}}$.
The number of steps over a sequence $(\e_k)\in\wt X_0^{\text{reg}}$
is equal to the length of the block beginning with $\e_1=1$.
\endproclaim
\demo{Proof} It suffices to present the steps of the corresponding
tower. Let, by definition,
$\wt X_1^{\text{reg}}=\wt\tau\wt X_0^{\text{reg}}$, and
$\wt X_{2j}^{\text{reg}}=\wt\tau \wt X_{2j-1}^{\text{reg}},\
\wt X_{2j+1}^{\text{reg}}=\wt\tau \wt X_{2j}^{\text{reg}}\setminus
(\e_1=1),\ j\ge1$. This completes the proof, as $\wt X^{\text{reg}}=
\bigcup_0^\infty \wt X_j^{\text{reg}}$, the union being disjoint.
\enddemo

\bigskip
\epsfysize=4cm
\centerline{\epsfbox{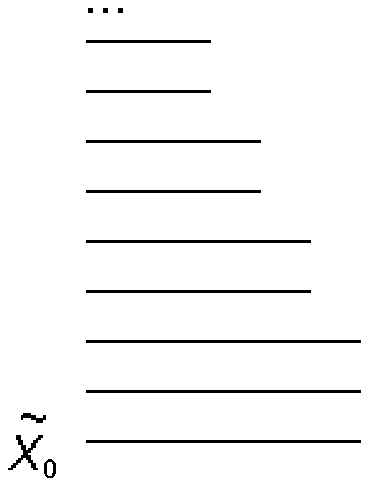}}
\bigskip\bigskip
\nopagebreak

\centerline{Fig. 2. The steps of the special automorphism $\wt\tau$}
\medskip

Below the corresponding result will be established for the metric case
with the two-sided Erd\"os measure.

\subhead 2.4. Bernoullicity of the goldenshift \endsubhead
In this subsection we will show that the goldenshift (one-sided
or two-sided) is a Bernoulli shift in the space $\X$ (respectively
$\wt\X$) with respect both to Lebesgue and Erd\"os measures
and compute their one-dimensional distributions.

We recall that
above we denoted the set of blocks by $\B$; let the same
letter stand for the totality of cylinder sets $\{B_1=B\}\subset \X$
for all blocks $B$. Let $\B_n$
denote all cylinder sets $\{B_1=B\}\subset\X$ such that
the length of $B$ is $2n+1$. Thus, $\B=\bigcup_{n=1}^\infty \B_n$.

Recall that the mapping $\Psi$ assigns to each regular $x$ its block
expansion, i.e. a certain sequence in the space $\X$.
We denote by  $m$ the normalized Lebesgue measure on the interval
$(\la^{-1},1)$; let $m_\X$ stand for the measure
$\Psi(m)$ in the space $\X$, and,  similarly, let $\mu_\X$ denote
$\Psi(\mu_{I_0})$ (this measure is well defined, as $\mu$-a.e. point
$x\in(\la^{-1},1)$ is also regular).

We recall that to any block we associated
the interval $\De_r$ defined as the image of the cylinder
``$B(r)1$'' in $X$ by the mapping~(1.1) (see Remark~1
after the definition of the goldenshift).

We are going to show that the measures $m_\X$ and $\mu_\X$
are Bernoulli in the space $\X$ and to compute their one-dimensional
distributions.
\proclaim{Theorem 2.4} The measure $m_\X$ in the space $\X$
is a product measure with equal multipliers, i.e. a Bernoulli measure.
\endproclaim
\demo{Proof} By the linearity of
$\S$ on each interval $\De_r$ and the fact that
$\S(\De_r)=(\la^{-1},1)$, we have for any Borel set $E\subset (\la^{-1},1)$,
$$
m(\S^{-1}E\cap\De_r)=mE\cdot m\De_r, \tag 2.1
$$
whence the required assertion immediately follows by virtue of the
obvious $\S$-invariance of $m$, and by setting $E=\De_{r'}$ in relation~(2.1)
for any $r'$, which yields the $m_\X$-independence of the first and
the second blocks.
\enddemo
So, it remains to compute the one-dimensional distribution of $m_\X$.

\proclaim{Corollary 2.5} For any cylinder set $\{B_1=B\}\subset \B_n$
its measure $m_\X$ equals $\la^{-2n-1}$.
The measure $m_\X$ of $\B_n$
is equal to $\frac1{2\la}\(\frac2{\la^2}\)^n$.
\endproclaim
\proclaim{Proposition 2.6} The measure $\mu_\X$ is a
product measure on $\X$ with equal multipliers.
\endproclaim
\demo{Proof} It suffices to establish a relation similar to (2.1) for
the measure $\mu_\X$ and for any finite block sequence $E=B_1\dots B_k$.
Note first that by virtue of Lemma~2.2, $\E(B_1\dots B_k)=\E(B_1)\dots
\E(B_k)$ for any blocks $B_1,\dots, B_k$.
Next,
$$
\n^{-1}(B_1\dots B_k1)=\E(B_1)\dots \E(B_k)\n^{-1}(\e_1=1).
$$
We are going to show that
$$
\mu_\X\{B_1\dots B_k\}=\mu_\X\{B_1\}\dots \mu_\X\{B_k\}
=\frac{f(B_1)}{2^{|B_1|}}\cdots \frac{f(B_k)}{2^{|B_k|}},\quad k\ge1.\tag 2.2
$$
To do this, we use previous remarks and the definition of the
Erd\"os measure on $X$ by means of the normalization (see Section~1).
We have $\mu_\X\{B_1\dots B_k\}=\mu(B_1\dots B_k1)/\allowmathbreak
\mu(\e_1=1)$,
and
$$
\align
\mu(B_1\dots B_k1)&=p(\n^{-1}(B_1\dots B_k1))=p(\E(B_1)\dots
\E(B_k)\n^{-1}(\e_1=1)) \\
                  &=\frac{f(B_1)}{2^{|B_1|}}\cdots \frac{f(B_k)}{2^{|B_k|}}
                  \cdot \mu(\e_1=1)
\endalign
$$
(by Lemma 2.2), whence the required assertion follows.
\enddemo

Thus, we have proved one of the main results of the present paper.
\proclaim{Theorem 2.7} The goldenshift $\S$ is a Bernoulli automorphism
with respect to Lebesgue and Erd\"os measures.
\endproclaim
Now we are ready to give the second proof of Erd\"os theorem
(see Section~1).
\proclaim{Corollary 2.8} {\rm (a new proof of Erd\"os theorem)}
The Erd\"os measure
is singular with respect to Lebesgue measure.
\endproclaim
\demo{Proof} In fact, by the corollary of the ergodic theorem,
the measures involved are
mutually singular on the interval $(\la^{-1},1)$ (recall
that $\mu(0,\la^{-2})=\frac13$, whence the measure $\mu$
differs from the Lebesgue measure).
This yields the assertion
of the corollary, as any infinite convolution of discrete measures,
if it is not stochastically constant,
is known to be either singular or absolutely continuous with respect
to Lebesgue measure (the ``Law of Pure Types'', see \cite{JeWi})).
\enddemo
The one-dimensional distribution of $\mu_\X$ is a bit more sophisticated
than for $m_\X$. It is described as follows (see formula (2.2)):

\proclaim{Proposition 2.9} For a block $B=B(a_1,\dots,a_t)$
of length $2n+1$,
$$
\mu_\X \{B_1=B\}=\frac{f(B)}{2^{|B|}}=\frac{p+q}{2^{2n+1}},
$$
where, as usual, $p/q=[a_1,\dots,a_t]$. The measure $\mu_\X$ of the set $\B_n$
equals $\frac13\cdot \(\frac34\)^n$.
\endproclaim

\subhead 2.5. Concluding remarks on the Erd\"os measure \endsubhead
We conclude the study of ergodic properties of Erd\"os measures
(one-sided and two-sided) and the transformations of shift and
goldenshift.

 Recall that the Erd\"os measure $\mu$ is quasi-invariant
under the one-sided shift $\tau$, and the equivalent measure $\nu$
is $\tau$-invariant. It is worthwhile to know the behavior
of $\nu$ with respect to the goldenshift $\S$.

Let $\nu_\X$ be defined in the same way as $\mu_\X$.
We formulate the following claim (for more details
see Appendix~C).
\proclaim{Proposition 2.10} The measure $\nu_\X$ on the space $\X$ of block
sequences is quasi-invariant under the goldenshift $\S$. More
precisely, any two cylinders $\{B_j=B'_j\}$ and $\{B_i=B'_i\}$ with
$i\neq j$ are $\nu_\X$-independent, and for $B=B(a_1,\dots,a_t)$,
$$
\align
\nu_\X \{B_k=B\}&=\mu_\X\{B_k=B\}=\frac{f(B)}{2^{|B|}}=
\frac{p+q}{2\cdot 4^{a_1+\cdots+a_t}},\quad k\ge2, \\
\nu_\X \{B_1=B\}&=
  \cases \frac{\frac45 p+\frac65 q}{2\cdot 4^{a_1+\cdots+a_t}},& B=100\dots \\
         \frac{\frac65 p+\frac45 q}{2\cdot 4^{a_1+\cdots+a_t}},& B=101\dots
  \endcases
\endalign
$$
\endproclaim
Now we will prove a numerical claim useful for the next section.
\proclaim{Corollary 2.11} $\wt\nu\wt X_0=\frac19$.
\endproclaim
\demo{Proof} We have by the definition of the set $\wt X_0$,
Proposition~2.10 and the fact that $\nu(\e_1=1)=
\frac5{18}$ (see Corollary~1.22),
$$
\wt\nu\wt X_0=\sum_{k=1}^\infty \wt\nu\bigl(1(00)^k1\bigr)=\frac5{18}
\sum_{k=1}^\infty \frac{\frac45+\frac65k}{1+k}\cdot\frac{k+1}{2^{2k+1}}=
\frac19.
$$
\enddemo
Finally, we prove a metric version of Proposition~2.3.

\proclaim{Theorem 2.12} The two-sided shift $\wt\tau$ with the
measure $\wt\nu$ is a special automorphism over the goldenshift
$\wt\S$ with the measure $\wt\mu$ on the space $\wt X_0$.
The step function is defined as the length of the block
beginning with the first coordinate.
\endproclaim
\demo{Proof} It suffices to show that the lifting measure for
$\wt\mu$ coincides with $\wt\nu$. This in turn is implied again by the
corollary of the
ergodic theorem applied to $\wt\tau$ that preserves both measures
which are ergodic.
Since they are clearly equivalent, we are done.
\enddemo

\head 3. The entropy of the goldenshift and applications \endhead

In this section we will establish a relationship between
the entropy of the Erd\"os measure in the sense
of A.~Garsia and the entropy of the goldenshift with
respect to $\mu$, i.e. between two different entropies.
As an application, we will reprove the formula for Garsia's entropy
proved in \cite{AlZa}. Further, we use random walk theory to
compute the dimension of the Erd\"os
measure on the interval.

\subhead 3.1. Fibonacci graph, random walk on it
 and Garsia's entropy \endsubhead
The combinatorics of equivalent 0-1 sequences may be expressed graphically,
namely by means of the Fibonacci graph introduced in \cite{AlZa}.
Let, as in Section~1, $\Sigma=\pr$, and
let the mapping $\pi:\Sigma\to[0,1]$ be defined as
$$
\pi(\e_1,\e_2,\dots)=\sum_{k=1}^\infty \e_k\la^{-k-1}.\tag 3.1
$$
Since the $\e_n$ assume the values 0 and 1 without any
restrictions, a typical $x$ will have a continuum number
of representations,
and they all may be illustrated with the help
of the {\it Fibonacci graph} depicted
in Fig.~3. This figure appeared for the first time in the work  due to
J.~C.~Alexander and D.~Zagier \cite{AlZa}.
Let us give the precise definition.

\definition{Definition} The Fibonacci graph $\Phi$ is a
binary graph with the edges labeled with 0 if to the left and 1 if to the right.
Any vertex at the $n$'th
level corresponds to a certain $x$, for which some representation~(3.1)
is finite with the length~$n$ (obviously, in this case $x=\{N\la\}$ for some
$N\in\Bbb Z$). The paths are 0-1 sequences treated as representations
of the form~(3.1).\footnote{The term ``Fibonacci graph'' is overloaded,
as the authors know several different graphs also called ``Fibonacci''.
Nevertheless, we hope that there will be no confusion with any of them.}

\bigskip
\epsfysize=6cm
\centerline{\epsfbox{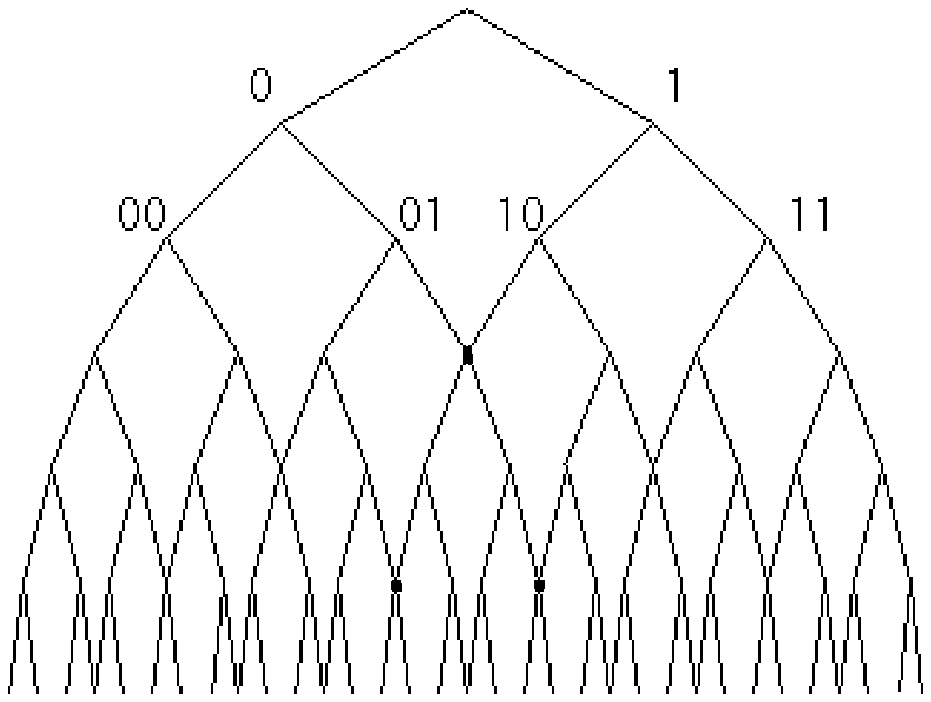}}
\bigskip\bigskip
\nopagebreak

\centerline{Fig. 3. The Fibonacci graph $\Phi$}
\medskip
\enddefinition

\remark{Remark} The vertices of the $n$'th level of the graph $\Phi$ can be
treated as the nonnegative integers from 0 to $F_{n+2}-2$. Namely,
if a path $(\e_1 \e_2\dots \e_n)$
goes to a vertex $k$, then, by definition,
$k=\sum_1^n \e_jF_{n-j}$ (obviously, this sum does not depend on the
choice of a path).
\endremark

Let $Y(\Phi)$ denote the set of paths in the graph $\Phi$. Obviously,
$Y(\Phi)$ is naturally isomorphic to $\Sigma$,
and sometimes we will not make a distinction between them. Let
$(\e_1\e_2\dots)$ be a path, and let the projection from
$Y(\Phi)$ onto $[0,1]$ be also denoted by $\pi$ (see formula~(3.1)).

Let, as above,  $f_n(k)$ denote the number of representations of a
nonnegative integer $k$ as a sum of not more than $n$ first Fibonacci
numbers. It is easy to see that $f_n(k)$ is also the {\it frequency}
of the vertex $k$ on the $n$'th level of the graph $\Phi$.\footnote{
Thus, relation~(1.2) completely determines the whole graph $\Phi$.}
Let $D_n=\{k: k=\sum_{k=1}^n \e_kF_{n-k},\ \e_k\in\{0,1\}\}$
(or, equivalently, the $n$'th level of the Fibonacci graph), and
$D'_n=\{w: w=\sum_{k=1}^n \e_k\la^{-k-1},\ \e_k\in\{0,1\}\}$. These sets
are clearly isomorphic ($w\lr k$), and
$\#D_n=\#D'_n=F_{n+2}-1$. The use of $D_n$ instead of
$D'_n$ is due only to technical reasons. We recall that the sequence
of distributions $(2^{-n}f_n(w))_{n=1}^\infty$
tends to the distribution of Erd\"os measure (see Section~1).

We have proved in the previous section that block is
an object defined on the space of almost all 0-1 sequences (not only
admissible), i.e. on the Fibonacci graph.
Let us now explain the geometric treatment of the block expansion.

Note first that
an odd level $2n+1$ contains $2^{n-1}$ specific vertices
which we will call, following \cite{AlZa}, the {\it Euclidean}
vertices. They are defined recursively. The central vertex on the
third level is Euclidean; then, any Euclidean vertex generates exactly
two new Euclidean vertices on the next odd level by means of the
arcs 00 and 11. We call the set of all Euclidean
vertices the {\it Euclidean tree}.

It is evident that for a given path in the
Fibonacci graph its first block can end only at some Euclidean vertex.
To find out, if it does end at a vertex $v$, we consider the induced Fibonacci
graph with $v$ as the top and also the induced Euclidean tree. The criterion
in question is that a given path should go through some induced Euclidean
vertex earlier than through any initial one. Next, considering the
induced Fibonacci graph, we can find the second block of the
block expansion, etc.

We denote the entropy of the discrete distribution on $D'_n$ decribed above,
by $H^{(n)}$. Thus,
$$
H^{(n)}=-\sum_{k=0}^{F_{n+2}-2}\frac{f_n(k)}{2^n}\log \frac{f_n(k)}{2^n}.
$$
Then, by definition,
$$
H_\mu:=\lim_{n\to\infty}\frac{H^{(n)}}{n\log\la}
$$
(this limit is known to exist and is obviously
independent of  the choice of the base
of logarithms, see \cite{Ga}). Now we choose
once and for all $\la$ as the
base of logarithms.

The quantity $H_\mu$ can be considered as the entropy of the
random walk on the Fibonacci graph with the probabilities
$(\frac12, \frac12)$.
In the next item it will be shown that in fact $H_\mu$ is
proportional to the entropy of the goldenshift.
The Erd\"os measure is the projection of the
Markov measure $(\frac12, \frac12)$
on the graph $\Phi$ under the mapping $\pi$. We consider
the random walk on the Fibonacci graph with the equal transition measures.

\definition{Definition} The {\it Fibonacci semigroup} (resp. {\it group})
is by definition, the semigroup (resp. group)
with the generators $a,b$ and the relation $ab^2=ba^2$.
\enddefinition
The following claims are straightforward.

\proclaim{Proposition 3.1} The Fibonacci graph is the Cayley graph
of the Fibonacci semigroup.
\endproclaim
In \cite{Av} (see also
\cite{KaVer}) was introduced the notion of the entropy of
a random walk on a finitely generated group (or semigroup).
The following claim establishes a relation between the two notions
of entropy. Note that the Fibonacci semigroup can be
naturally embedded into the Fibonacci group, that is why
we can use the theory of random walks on groups.

\proclaim{Proposition 3.2} The entropy $H_\mu$ is equal to
the entropy of the random walk on the Fibonacci semigroup with
the probabilities $\(\frac12,\frac12\)$.
\endproclaim
\demo{Proof} Follows from the definitions and Proposition~3.1.
\enddemo

\subhead 3.2. Main theorem \endsubhead
We prove an assertion that is one of the central points of the present paper.
Note that in \cite{AlZa} Garsia's entropy was computed by means of
generating functions. We will see that $H_\mu$ is closely
connected with the entropy of the goldenshift, which gives a new
and simplified proof of their relation and relates it to the dynamics
of the Erd\"os measure.

\proclaim{Theorem 3.3} The following relation holds:
$$
h_\mu(\S)=9H_\mu.
$$
\endproclaim
\demo{Proof} We are going to apply Abramov's formula
for the entropy of the special automorphism (see \cite{Ab})
to the dynamical systems $(\wt X,\wt\nu,\wt\tau)$ and
$(\wt X_0,\wt\mu,\wt\S)$. By Abramov's formula and Theorem~2.12,
$$
h_{\wt\mu}(\wt\S)=\frac1{\wt\nu\wt X_0} h_{\wt\nu}(\wt\tau).
$$
From this relation we will deduce the required one.\newline
1. Since the dynamical system $(\wt X, \wt\mu, \wt\S)$ is the natural
extension of $(X, \mu, \S)$, we have $h_{\wt\mu}(\wt\S)=h_\mu(\S)$
and by the same reason,
$h_{\wt\nu}(\wt\tau)=h_\nu(\tau)$. \newline
2. By Corollary~2.11, $\wt\nu\wt X_0=\frac19$.\newline
3. The rest of the proof is devoted to establishing the validity
of the relation
$$
h_\nu(\tau)=H_\mu.
$$
 Let $\eta$ denote
the partition of $X$ into the cylinders $(\e_1=0)$ and $(\e_1=1)$.
Since $\eta$ is a generating partition for $\wt\tau$, we
have $h_\nu(\tau)=h_\nu(\tau,\eta)$ by
Kolmogorov's theorem. By definition,
$$
h_\nu(\tau,\eta)=\lim_{n\to\infty}\frac 1n H_\nu(\eta^{(n)}),
$$
where $\eta^{(n)}$ is the partition of $X$ into $F_{n+1}$ admissible
cylinders of the form $(\e_1=i_1,\dots, \allowmathbreak\e_n=i_n)$.
We need to prove that
$$
H_\nu(\eta^{(n)})\sim H^{(n)}.
$$
By virtue of the equivalence of
the measures $\mu$ and $\nu$ it suffices to show this for
$H_\mu(\eta^{(n)})$ instead of $H_\nu(\eta^{(n)})$.
Let $\theta_n(k)=2^{-n}f_n(k)$. Then, by relation~(1.2),
for $n\ge 3$,
$$
\theta_n(k)= \cases \frac12\theta_{n-1}(k), & 0\le k\le F_n-1 \\
               \frac12(\theta_{n-1}(k) + \theta_{n-1}(k-F_n)),
               & F_n\le k\le F_{n+1}-2 \\
               \frac12\theta_{n-1}(k-F_n), & F_{n+1}-1 \le k \le F_{n+2}-2.
        \endcases
$$
We are going to obtain almost the same
recurrence relation for the distribution~$\eta^{(n)}$. To do this,
we return to the interval $[0,1]$ and denote by $\eta^{(n)}$ the
partiton into $F_{n+1}$ intervals which is the image of the corresponding
partition of $X$ with the help of the  canonical expansion.
So, let $\eta^{(n+1)}=:(J_n(k))_{k=0}^{F_{n+2}-1}$ with the ordered
intervals $J_n(k)$. Finally, let $\mu_n(k):=\mu J_n(k)$. Then
by Lemma~1.1, for
$n\ge 3$,
$$
\mu_n(k)= \cases \frac12\mu_{n-1}(k), & 0\le k\le F_n-1 \\
               \frac12(\mu_{n-1}(k) + \mu_{n-1}(k-F_n)),
               & F_n\le k\le F_{n+1}-1 \\
               \frac12\mu_{n-1}(k-F_n), & F_{n+1} \le k \le F_{n+2}-3.
        \endcases
$$
Also, $\mu_n(F_{n+2}-2)=O(2^{-n}),\,  \mu_n(F_{n+2}-1)=O(2^{-n})$.
Thus, by induction on $n$, there exists $C>0$ such that
$$
\frac 1C\le \frac{\theta_n(k)}{\mu_n(k)}\le C, \quad 0\le k\le F_{n+2}-2,
$$
Since the measure $\nu$ is $\tau$-invariant and ergodic,
we can apply the Shannon-McMillan-Breiman theorem and deduce that the
entropies of the distributions $\mu_n(k)$ and $\theta_n(k)$ are
equivalent. \qed
\enddemo

\subhead 3.3. Alexander-Zagier's theorem \endsubhead
We first compute the entropy of the goldenshift with
respect to the Erd\"os and Lebesgue measures.

\proclaim{Proposition 3.4} The metric entropies of the goldenshift $\S$ with
respect to the two measures  in question are computed as follows:
$$
\align
h_m(\S)   &=-\sum_{B\in \B}m_\X(B)\log_\la m_\X(B)=4\la+3=9.4721356\dots,  \\
h_\mu(\S) &=-\sum_{B\in \B}\mu_\X(B)\log_\la \mu_\X(B)=-
\sum_{n=1}^\infty\sum_{B:|B|=2n+1}\frac{p+q}{2^{2n+1}}\
\log_\la\frac{p+q}{2^{2n+1}}=
8.961417\dots
\endalign
$$
\endproclaim
\demo{Proof} This is a direct computation using
the Bernoullicity of $\S$ with respect to both measures
and Corollaries~2.5 and 2.9.
\enddemo

\proclaim{Corollary 3.5}
Let $\La:=\log_\la 2$, and
$$
k_n=\sum \Sb t\ge 1 \\ (a_1,\dots,a_t)\in\Bbb N^t \\ a_1+\cdots+a_t=n
         \\p/q=[a_1,\dots,a_t]
         \endSb (p+q)\log_\la(p+q).
$$
Then
$$
h_\mu(\S)=9\(\La-\frac1{18}\sum_{n=1}^\infty\frac{k_n}{4^n}\).\tag 3.2
$$
\endproclaim
\demo{Proof} We have
$$
\align
h_\mu(\S) &=-\sum_{n=1}^\infty\sum_{B:|B|=2n+1}\frac{p+q}{2^{2n+1}}\
\log_\la\frac{p+q}{2^{2n+1}} \\
          &=-\frac12\sum_{n=1}^\infty 4^{-n}\(k_n-\La(2n+1)
\sum_{a_1+\cdots+a_t=n}(p+q)\).
\endalign
$$
It suffices now to recall that
$$
\sum \Sb t\ge 1 \\ (a_1,\dots,a_t)\in\Bbb N^t \\ a_1+\cdots+a_t=n
         \\p/q=[a_1,\dots,a_t]
         \endSb (p+q)=2\cdot 3^{n-1}
$$
and to compute the value of the corresponding series.\qed
\enddemo
Note that
the quantity $k_n$ appeared for the first time in \cite{AlZa}
in somewhat different notation. Namely, let $k$ and $i$ be positive integers,
and let $e(k,i)$ denote the length of the simple Euclidean algorithm
for $k$ and $i$ (formally: $e(i,i)=0,\,e(i+k,i)=e(i+k,k)=e(i,k)+1$).
Then obviously
$$
k_n=\sum \Sb k>i>0 \\ \gcd(k,i)=1,\,\, e(k,i)=n
         \endSb
k\log_\la k.
$$
In the cited work J.~C.~Alexander and D.~Zagier used this
definition of $k_n$ to deduce a formula for $H_\mu$ in terms of $k_n$.
We will prove their assertion in two different ways. The first is an
immediate consequence of Theorem~3.3 and relation~(3.2),
while the second is rather long
but reveals a more essential relationship between certain structures on the
Fibonacci graph (see Appendix~D).

\proclaim{Proposition 3.6} {\rm (Alexander-Zagier, 1991)}.
The following relation holds:
$$
H_\mu=\La-\frac1{18}\sum_{n=1}^\infty\frac{k_n}{4^n}=0.995713\dots
\tag3.3
$$
\endproclaim
\demo{Proof}  An application of Theorem~3.3 and of
relation~(3.2).
\enddemo

\subhead 3.4. The dimension of the Erd\"os measure \endsubhead
As an application of the treatment of the Erd\"os measure
as the projection of the measure of the uniform random
walk on the Fibonacci graph, we will compute the dimension
of $\mu$ in the sense of L.-S.~Young.

We first give
a number of necessary definitions (see \cite{Y}).
\definition{Definition} Let $\nu$ be a Borel probability measure on
a compact space $Y$. The quantities
$$
\align
\dim_H \nu &=\inf\{\dim_H A: A\subset Y,\,\nu A=1\}, \\
\ov C(\nu) &=\limsup_{\de\to0}\,\inf\{\ov C(A): A\subset Y, \,
\nu A\ge 1-\de\},\\
\un C(\nu) &=\liminf_{\de\to0}\,\inf\{\un C(A): A\subset Y, \,
\nu A\ge 1-\de\}
\endalign
$$
(where $\ov C(A)$ and $\un C(A)$ are respectively the upper and lower
capacities of $A$) are called the {\it Hausdorff dimension} of a measure $\nu$
and
the {\it upper} and {\it lower capacities} of $\nu$, respectively.
\enddefinition
Let next $N(\e,\de)$ denote the minimal number of balls of radius $\e>0$
which are necessary to cover a set of $\nu$-measure $\ge 1-\de$.
\definition{Definition} The quantities
$$
\align
\un C_L(\nu)&=\limsup_{\de\to0}\,\,\liminf_{\e\to0}\frac{\log N(\e,\de)}
{\log (1/\e)}, \\
\ov C_L(\nu)&=\limsup_{\de\to0}\,\,\limsup_{\e\to0}\frac{\log N(\e,\de)}
{\log (1/\e)}
\endalign
$$
are called the {\it lower} and {\it upper Ledrappier capacities} of $\nu$.
\enddefinition
\definition{Definition} Let $H_\nu(\e)=\inf\{H_\nu(\xi):\diam\xi\le \e\}$,
where $H_\nu(\xi)$ is the entropy of a finite partition $\xi$. The quantities
$$
\align
\ov R(\nu)&=\limsup_{\e\to0}\frac{H_\nu(\e)}{\log(1/\e)},\\
\un R(\nu)&=\liminf_{\e\to0}\frac{H_\nu(\e)}{\log(1/\e)}
\endalign
$$
are called respectively the {\it upper} and {\it lower informational
dimensions} of $\nu$ ($=$ {\it R\'enyi dimensions}).
\enddefinition

\proclaim{Theorem} {\rm (L.-S. Young \cite{Y}, 1982).} Let $\nu$
be a Borel probability measure on a metric space $Y$, and let
$B(x,r)$ denote the ball with the center at $x$ of radius $r$. If
$$
\a(x):=\lim_{r\to0}\frac{\log\,\nu B(x,r)}{\log r}\equiv \a
$$
for $\nu$-a.e. point $x\in Y$, then
$$
\dim_H\nu=\ov C(\nu)=\un C(\nu)=\ov C_L(\nu)=\un C_L(\nu)=\ov R(\nu)=
\un R(\nu)=\a.
$$
\endproclaim
\definition{Definition} If the condition of Young's theorem
is satisfied, then this $\a$ is called the {\it pointwise dimension} of
a measure $\nu$ and is denoted by $\dim(\nu)$. This
notion was also proposed in \cite{Y}.
\enddefinition

\proclaim{Theorem 3.7} For the Erd\"os measure $\mu$,
$$
\dim(\mu)=H_\mu.
$$
\endproclaim
\demo{Proof}
Fix a path $\ov\e=(\e_1\e_2\dots)\in Y(\Phi)$, and let
$x=\pi(\ov\e)$ and $Y_n=Y_n(\ov\e)$ be the interval whose
every point $x'$ has a path ${\ov\e}'\in Y(\Phi)$ such that
$\e_i\equiv\e'_i,\ 1\le i\le n$. Clearly,
$Y_n(\ov\e)=\[\sum_1^n\e_k\la^{-k-1},\sum_1^n\e_k\la^{-k-1}+\la^{-n}\]$,
We recall that $\mu=\vartheta_1*\vartheta_2*\dots$ (see Section~1).
Let $\mu^{(n)}:=\vartheta_1*\dots*\vartheta_n$.
Then by Shannon's theorem for the random walks (see
Theorem~2.1 in \cite{KaVer} and also \cite{De})
and because $H_\mu$ is the entropy of the random walk on the
Fibonacci semigroup (see Proposition~3.2),
$$
\lim_{n\to\infty}\frac{\log_\la\mu^{(n)}(Y_n(\ov\e))}n =-H_\mu
$$
for $\mu$-a.e. $x\in[0,1]$. We set
$Y_n'(\ov\e)=\[\sum_1^n\e_k\la^{-k-1},\sum_1^n\e_k\la^{-k-1}+\la^{-n+2}\]
\supset Y_n(\ov\e)$. Obviously,
$\mu^{(n)}(Y_n'(\ov\e))\sim\mu^{(n)}(Y_n(\ov\e))$.
Given $h>0$, we choose $n=n(h)$ such that
$Y_{n+1}'(\ov\e)\subset (x, x+h)\subset Y_n'(\ov\e)$
for any $\ov\e\in\pi^{-1}\{x\}$. Hence it follows that for
$\mu$-a.e. $x$,
$$
\lim_{h\to0}\frac{\log_\la\mu(x,x+h)}{\log_\la h}=
-\lim_{n\to\infty}\frac{\log_\la\mu Y_n(\ov\e)}n=H_\mu,
$$
as $h\asymp \la^{-n}$. \qed
\enddemo
\remark{Remark} When the present paper was in preparation,
the authors were told that the claim of Theorem~3.7 can be
obtained as a corollary of several results including the new one
due to F.~Ledrappier and A.~Porzio.
More precisely, it was shown in \cite{AlYo} that
$H_\mu=\ov R(\mu)=\un R(\mu)$, and in \cite{LePo} it was
proved that the limit in Young's theorem does exist for the Erd\"os measure.
This proves Theorem~3.7.

Our proof is straightforward and, what is more important, is
a direct corollary of a Shannon-like theorem, so far it leads to
new connections between geometric and dynamical properties of the Erd\"os
measure.
\endremark

\proclaim{Corollary 3.8}
$$
\dim_H\mu=\ov C(\mu)=\un C(\mu)=\ov C_L(\mu)=\un C_L(\mu)=\ov R(\mu)=
\un R(\mu)=H_\mu=0.995713\dots
$$
\endproclaim
\remark{Remark {\rm1}} Another proof of Theorem~3.7 can be obtained
by using the Bernoulli structure of the measure $\mu$.
More precisely, for a regular $x\in(\la^{-1},1)$ having a {\bf normal}
block expansion with respect to the measure $\mu_\X$, as is easy
to compute, the limit in the definition of the dimension
equals $\frac19 h_\mu(\S)$. This yields also another proof of
Theorem~3.3. The details are left to the interested reader.
\endremark
\remark{Remark {\rm2}} In fact, we have computed the $\mu$-typical Lipschitz
exponent of the distribution function of the Erd\"os measure.
Note that in \cite{Si} it was proved that the best possible Lipschitz
exponent of this function for all $x$ is $\La-\frac12=0.9404\dots$.
\endremark
\remark{Remark {\rm3}} As a conjecture, we claim that for the two-dimensional
Erd\"os measure $\wt\mu$ (see Section~1),
$$
\dim(\wt\mu)=2\dim(\mu)=1.991426\dots
$$
The proof could apparently follow from the theorem due to L.-S.~Young
\cite{Y} relating the pointwise dimension to the entropy of an automorphism
and Lyapunov exponents of an ergodic measure. Besides, we think that
the measure $\wt\mu$ has a local structure of direct product, which
would also explain this relation.
\endremark

\subhead 3.5. Fibonacci exponents \endsubhead
We recall that in Section~2 we defined the function $f$ acting on the set
of finite 0-1 words and counting the number of words equivalent
to an argument. Let $x\in(0,1)$ and $(\e_1\e_2\dots)$ be its
canonical representation. Let now the finite word $x_n=(\e_1\dots\e_n)$.
\definition{Definition} The limit
$$
E(x)=\lim_{n\to\infty}\frac{\log_\la f(x_n)}n
$$
(if exists) will be called the {\it Fibonacci exponent} of a point $x$.
\enddefinition
We will show that for a.e. $x$ with respect to Erd\"os measure the
Fibonacci exponent exists and is the same, as well as for Lebesgue measure.
Besides, we reprove one theorem due to S.~Lalley (see \cite{Lal} and
references therein). The proof in \cite{Lal} can be applied to any PV
number $\la$ but for the golden ratio our proof is more direct.

\proclaim{Proposition 3.9} {\rm (\cite{Lal, Th. 2}).} For $\mu$-a.e. $x$,
$$
E(x)= E_\mu := \La-H_\mu=0.44469\dots
$$
\endproclaim
\demo{Proof} It suffices to consider a regular $x\in (\la^{-1},1)$.
Let $B_1B_2\dots$ be its block expansion.
Then
$$
E(x)=\lim_{k\to\infty} \frac{\log_\la f(B_1\dots B_k)}{|B_1|+\cdots+|B_k|}=
\frac{\EE_\mu \log_\la f(B_1)}{\EE_\mu |B_1|}
$$
(where $\EE_\mu$ denotes mathematical expectation with respect to Erd\"os
measure)
by the ergodic theorem applied to the goldenshift and Erd\"os measure.
Now it suffices to observe that $\EE_\mu |B_1|=9$, and
$\EE_\mu \log_\la f(B_1)=\frac12\sum_{n\ge1}\frac{k_n}{4^n}$
and apply Proposition~3.6.\qed
\enddemo
In the same way we obtain

\proclaim{Proposition 3.10} For a.e. $x$ with respect to Lebesgue measure,
$$
E(x):=E_m = \frac{\EE_m\log_\la f(B_1)}{\EE_m |B_1|}=
\frac{\sum_{n=1}^\infty \ell_n \la^{-2n-1}}{4\la+3},
$$
where
$$
\ell_n=\sum \Sb t\ge 1 \\ (a_1,\dots,a_t)\in\Bbb N^t \\ a_1+\cdots+a_t=n
         \\p/q=[a_1,\dots,a_t]
         \endSb \log_\la(p+q).
$$
\endproclaim
\remark{Remark {\rm1}} We established in Proposition~3.9 a relation between
the pointwise dimension of the Erd\"os measure and its typical Fibonacci
exponent. This gives us an occasion to state without
proof a similar claim for Lebesgue
measure:
$$
\lim_{h\to0} \frac{\log \mu(x,x+h)}{\log h}=\La - E_m
$$
for a.e. $x$ with respect to the Lebesgue measure. The proof is the same
as the one mentioned in Remark~1 after Corollary~3.8. Apparently, this
dimensional characteristic has not been considered yet.
\endremark
\remark{Remark {\rm2}} Let $f(n)$ be the number of representations of
a positive integer $n$ as a sum of distinct Fibonacci numbers.
If $n=\sum_1^k \e_j F_j$ is such a representation with
$\e_j\in\{0,1\}\ \e_j\e_{j+1}=0,\ 1\le j\le k-1$, and
$\e_k=1$, then evidently $f(n)=f(\e_k\dots\e_1)$ in the
usual sense, which explains the choice of the
notation. It is known that $f(n)=O(\sqrt n)$ is
an attainable estimate (see \cite{Pu}) and that
the average behavior of $f$ in the sense of its summation function is
$\sum_{n<N} f(n) \asymp N^\La$ (for more precise results see \cite{DuSiTh}).
It is worth asking the question about the ``typical exponent'' of $f(n)$
with respect to density, i.e.
$$
E_d=\lim \Sb n\to\infty \\ n\in J \endSb\frac{\log f(n)}{\log n},
$$
where $J\subset \Bbb N$ is a subsequence of
density 1. We conjecture that this exponent
exists, and $E_d=E_m$,
i.e. density 1 corresponds to full Lebesgue measure.
\endremark

\head Appendix A.
The ergodic central measures and the adic transformation
on the Fibonacci graph \endhead

In this appendix we will study in detail some properties of
the space of paths $Y(\Phi)$ of the Fibonacci
graph introduced in Section~3. We first give
a necessary definition which is close to the defintion
of canonical expansions but reflects the fact that
 0 and 1 have the same rights in the graph $\Phi$.
\definition{Definition} The {\it generalized} canonical
expansion of a point in $(0,1)$ is defined as follows.
We construct the sequence $(\e_1 \e_2\dots)$ such that
relation~(3.1) holds, and either $(\e_1\e_2\dots)\in X$,
or $\e_1=\dots=\e_m=1$ for some $m\in \Bbb N$, and
the tail is in $X$. The algorithm is a clear modification
of the greedy algorithm.
\enddefinition

\remark{Remark} In fact, the generalized canonical expansions
lead to a normal form in the semigroup corresponding
to the group $G$ (see Section~3).
\endremark

The {\it tail partition} $\eta(\Phi)$ of $Y(\Phi)$ is
defined as follows.
\definition{Definition} Paths
$(\e_n)$ and $(\e_n')$, by definition,
belong to one and the same element of $\eta(\Phi)$ iff
\roster
\item "(i)" $\pi(\e_1\e_2\dots)=\pi(\e'_1\e'_2\dots)$, and
\item "(ii)" there exists
$N\in\Bbb N$ such that $\e_n\equiv\e_n',\, n>N$.
\endroster
\enddefinition
The partial lexicographic ordering on $Y(\Phi)$ is
defined for paths belonging to one and the same element of the tail
partition $\eta(\Phi)$.
\definition{Definition} Let two paths $\ov\e=(\e_1\e_2\dots)$
and $\ov\e'=(\e'_1\e'_2\dots)$ belong to one and the same element
of $\eta(\Phi)$. If $\e_{k-1}=0,\,\e_k=1,\,\e_{k+1}=1$ and
$\e'_{k-1}=1,\,\e'_k=0,\,\e'_{k+1}=0$ for some $k\ge2$, and
$\e_j\equiv\e'_j$ for $k-j\ge 2$, then, by definition,
$\ov\e\prec \ov\e'$. Next, by transitivity, $\ov\e\prec \ov\e',\
\ov\e'\prec \ov\e''$ implies $\ov\e\prec \ov\e''$.

\remark{Remark} This definition is consistent, because any element of
$\eta(\Phi)$ is isomorphic to a finite number of finite paths,
and they all can be transfered into one another with the help of
replacements $011\leftrightarrow 100$. Note also that this linear
ordering on each element of $\eta(\Phi)$ is stronger than the
partial ordering introduced in \cite{Pu} (see the end of item~2.2).
For example, $(100011*)\prec (011100*)$ in the above sense but in the
sense of the partial order they are noncomparable.
\endremark
\enddefinition
\definition{Definition} The {\it adic} transformation $T_\Phi$
assigns (if possible) to a path $\ov\e\in Y(\Phi)$ the path
$\ov\e'$ such that $\ov\e'$ belongs to the same
element of the tail partition as $\ov\e$ and is the immediate successor of
$\ov\e$ in the sense of the lexicographic order.
\enddefinition

It is clear that the adic transformation $T_\Phi$
is  not everywhere well defined.
More precisely, it is well defined on the paths $\ov\e$ containing
at least one triple $\e_{k}=0,\, \e_{k+1}=1,\, \e_{k+2}=1$.
Let us describe its action in more detail.
Let $(\e_1\e_2\dots)\in Y(\Phi)$ be as described. After finding
the first triple $\e_{k}=0,\, \e_{k+1}=1,\, \e_{k+2}=1$,
we \newline
1) replace it by $\e_{k}=1,\, \e_{k+1}=0,\, \e_{k+2}=0$, \newline
2) leave the tail $(\e_{k+3}\e_{k+4},\dots)$ without changes, \newline
3) find the minimal possible $(\e'_{1}\dots \e'_{k-1})$
equivalent to $(\e_1\dots \e_{k-1})$ in the sense of Section~2.

To carry out 3), we may use the algorithm of ``anti-normalization'',
i.e. the process analogous to the ordinary normalization but
changing ``100'' to ``011'' (cf. Section~1).

So, the generalized canonical expansions are just the {\it maximal} paths,
i.e. the ones where $T_\Phi$ is not well defined; thus, the set of maximal
paths is naturally isomorphic mod 0 to the interval $[0,1]$.
Geometrically the generalized canonical expansion corresponds
to the right most possible path
descending to a given vertex. Similarly,
the {\it minimal} paths (i.e. the ones, on which $T_\Phi^{-1}$
is not well defined) are just the so-called {\it lazy} expansions
(for the definition see, e.g., \cite{ErJoKo}).

For more general definitions of adic transformations and investigation
of their properties see \cite{Ver1}, \cite{Ver2}, \cite{LivVer} and
\cite{VerSi}.

Let us formulate several well-known definitions related to graded
graphs (see \cite{StVo} and \cite{VerKe} for more details).
We recall that topologically the space $Y(\Phi)$ is a nonstationary
Markov compactum (see \cite{Ver2} for definition).

\definition{Definition} A measure $\xi$ on $\Phi$ with the distribution
$\xi_n$ on its $n$'th level is called {\it Markov} if the sequence
$(D_n,\xi_n)$ of random variables is a (nonhomogenous) Markov chain.
\enddefinition
Now we can define for a Markov measure the notion of conditional measures.

\definition{Definition} A Markov measure $\xi$ on the graph $\Phi$ is called
{\it central} if any of the following equivalent conditions is satisfied:
\roster
\item For any vertex in this graph the conditional
measure on the set of all paths descending
to this vertex, is uniform.
\item $\xi$ is $T_\Phi$-invariant.
\item $\xi$ is invariant with respect to the tail partition $\eta(\Phi)$.
\endroster
\enddefinition
\definition{Definition} A central measure $\xi$
on $\Phi$ is called {\it ergodic}
if either of the following two equivalent conditions is satisfied:
\roster
\item The adic transformation $T_\Phi$ is ergodic with respect to it.
\item The tail partiton is $\xi$-trivial, i.e. contains only sets whose
$\xi$-measure is either 0 or 1.
\endroster
\enddefinition

The aim of this section is to describe \newline
1) all ergodic central measures on $\Phi$. \newline
2) the action of the adic transformation $T_\Phi$ on $\Phi$.

In the following
theorem we will describe the ergodic central measures
and the corresponding components of the action of $T_\Phi$.
As was noted above, $T_\Phi$ interchanges representations of one
and the same $x$.
We will see that the regularity or irregularity of
the generalized canonical expansion
of a given $x$ leads to three types of possible ergodic components
of the action of $T_\Phi$, namely, to a ``full'' odometer,
an irrational rotation of the circle or a special automorphism
over a rotation.

\proclaim{Theorem A.1} 1. The ergodic central measures on $\Phi$ are
naturally parametrized by the points of the interval $[0,1]$.
We denote by $\mu_x$ the measure corresponding to $x$. \newline
2. The measure $\mu_x$ is continuous if and only if $x\neq\{N\la\}$ for
any $N\in\Bbb Z$ (or, equivalently, if $x$ has infinite canonical
expansion). \newline
3. The action of the adic transformation $T_\Phi$ is not transitive, and its
trajectories are described as follows. Let $x$ be as in the
previous item, and let $\phi_x$ denote the space of paths in $\Phi$ such
that $\phi_x=\supp\mu_x$. The set $\phi_x$ is invariant
under $T_\Phi$ and we have the following alternatives.
\roster
\item "{\bf a.}" If the generalized
canonical expansion of $x$ contains infinitely
many pieces ``$1(00)^l1$'' with $l\ge 1$  (let us call
such a piece {\bf even}), then $T_\Phi|_{\phi_x}$ is strictly ergodic and
metrically isomorphic to the shift by $1$
on the group of certain {\bf a}-adic
integers with $\bold a=(a_1,a_2,\dots)$ being a sequence of positive
integers, generally speaking, nonstationary.
Thus, $T_\Phi|_{\phi_x}$ has a purely discrete rational spectrum.
\item "{\bf b.}" If the generalized
canonical expansion of $x$ does not contain
even pieces at all, then $T_\Phi|_{\phi_x}$ is also strictly ergodic and
metrically isomorphic to a certain irrational rotation of the circle.
\item "{\bf c.}" Finally, if the generalized
canonical expansion of $x$ contains
a finite number of even pieces,  then $T_\Phi|_{\phi_x}$ is metrically
isomorphic to some special automorphism over a rotation of the circle,
i.e. to a shift on the space $S^1\ltimes \Bbb Z/k$ for some $k\in\Bbb N$.
\endroster
\endproclaim
\demo{Proof} (1) Let $\phi_x$
be the set of all paths projecting to $x\in [0,1]$ (a $\pi$-fiber over $x$).
Obviously, the set $\phi_x$ is invariant under $T_\Phi$ for any $x$.
Thus, $T_\Phi$
is not transitive, and its action splits into components, each
acting on a certain $\phi_x$ (below we will see that for all $x$, except for
some
countable set, the action of $T_\Phi|_{\phi_x}$ is strictly
ergodic). \newline
(2) We have the following cases. If $x=0$ or $x=1$, then $\#\phi_x=1$.
If $x=\{N\la\}$ for some $N\in\Bbb Z$, then it is easy to see that
$\phi_x$ is countable and that the unique invariant measure for $T_\Phi$
is concentrated in a finite number of paths
(see Example~1 below). Henceforward in this proof
we assume that $x$ has infinite canonical expansion.
Let $x=\sum_{j=1}^\infty \e_j\la^{-j-1}$ be the generalized
canonical expansion of $x$.
We first split it in the following way: $(\e_1 \e_2 \dots)=B^{(0)}B^{(1)}$,
where $B^{(0)}$ is either $0^s$ or $1^s$ for some $s\ge 0$ (if $s=0$,
then $B^{(0)}=\emptyset$), and $B^{(1)}$ begins with ``$10$''.
Such a splitting is caused by the trivial reason that the action of $T_\Phi$ does
not touch at all the set $B^{(0)}$, as $T_\Phi$ only interchanges certain
triples ``$100$'' and ``$011$''.
So, we have the
following cases (they correspond to those enumerated in the theorem).
\newline
{\bf a.} If $B^{(1)}$ contains infinitely many even pieces, then
$x$ is regular (see Section~2), hence,
$$
B^{(1)}=B_1B_2B_3\dots,
$$
where
$$
B_j=1(00)^{a_1^{(j)}}(01)^{a_2^{(j)}}(00)^{a_3^{(j)}}
\dots(00)^{a_{t_j}^{(j)}}
\quad \text{or} \quad
B_j=1(01)^{a_1^{(j)}}(00)^{a_2^{(j)}}(01)^{a_3^{(j)}}
\dots(00)^{a_{t_j}^{(j)}}
$$
with $a_i^{(j)}\in \Bbb N$ and $t_j<\infty$.
\newline {\bf b.} If $B^{(1)}$ does not contain any even piece,
then obviously
$$
B^{(1)}=1(00)^{a_1}(01)^{a_2}(00)^{a_3}\dots \quad \text{or} \quad
B^{(1)}=1(01)^{a_1}(00)^{a_2}(01)^{a_3}\dots
$$
with $a_j\in\Bbb N$ for any $j\ge1$. \newline
{\bf c.} Finally, if the number of even pieces is finite (but nonzero),
then
$$
B^{(1)}=B_1B_2\dots B_m\wt{B}^{(1)},
$$
where $B_1,\dots, B_m$ have the form described in the previous item, and
$\wt{B}^{(1)}$ is an infinite block of the form described in item~{\bf b}.
\newline (3) Consider items~{\bf a}, {\bf b}, {\bf c}
from the viewpoint of the action of $T_x:=T_\Phi|_{\phi_x}$. \newline
{\bf a.} The idea of the study of $T_x$ in this case is based on two
assertions of the previous section, namely on Lemmas~2.1 and 2.2.
In particular, from Lemma~2.2 it follows
that blocks $B_i$ and $B_{i+1}$ for
all $i\in\Bbb N$ are replaced by any equivalent sequences independently, and
hence it is clear that for such a point $x$ the transformation~$T_x$
is the shift by 1 in the group of {\bf a}-adic integers with
$\bold a=(p_{1}+q_{1}, p_{2}+q_{2}, \dots$). This
transformation $T_x$ is
known to be strictly ergodic, i.e. there is a unique (product)
measure $\mu_x$ invariant under it. \newline
{\bf b.} We recall that in this case
$$
B^{(1)}=1(00)^{a_1}(01)^{a_2}(00)^{a_3}\dots \quad \text{or} \quad
B^{(1)}=1(01)^{a_1}(00)^{a_2}(01)^{a_3}\dots
$$
Let $\a=[1, a_1, a_2, \dots]$ denote a (regular) continued fraction.
We claim that in this case the transformation~$T_x$ is strictly
ergodic and metrically isomorphic to the rotation through~$\a$.
The idea of the proof lies in recoding the space $\phi_x$ into
the second model of the adic realization of the rotation from
\cite{VerSi} (see Section~3 of the cited work and Example~3 below).
The unique invariant measure can be described with the help
of Theorem~2.3 from the cited work. \newline
{\bf c.} This case in a sense is a ``mixture'' of the previous ones.
One can easily see that if $\wt{B}^{(1)}$ is parametrized by the
infinite sequence $(a_1, a_2, \dots)$ in the sense of the previous
item, and if $B_j=B_j(a_1^{(j)}, \dots, a_{t_j}^{(j)})$, then
$T_x$ acts on $\phi_x$ as the special automorphism over the
rotation through $\a=[1, a_1, a_2, \dots]$ with the constant
step function ($\equiv 1$) and the number of upper steps equal
to $\prod_1^m (p_{j}+q_{j})-1$. So, $T_x$ is again strictly ergodic.
The proof of the theorem is complete.
\enddemo
\remark{Remarks} 1. It is known (see, e.g., \cite{VerKe}) that
any ergodic central measure on the Pascal graph is also parametrized
by a real in $[0,1]$ but in a completely different manner, namely
by means of the first transition measure. It is appropriate to
compare that situation with the Fibonacci graph. We see that in the
graph~$\Phi$, for any $\a\in[0,1]$ there exists a central ergodic measure
$\mu$ such that $\mu(\e_1=0)=\a$. If $\a$ is irrational, then this
measure is unique, namely $\mu=\mu_x$ for $x=\sum_j\e_j\la^{-j-1}$
with $\a=[1, a_1, a_2,\dots]$ and
$(\e_1 \e_2 \dots)=1(00)^{a_1}(01)^{a_2}\dots$ for $\a>\frac12$, and
$1-\a=[1, a_1, a_2,\dots]$ and
$(\e_1 \e_2 \dots)=1(01)^{a_1}(00)^{a_2}\dots$ otherwise.
If $\a$ is rational, then there exists a whole interval of $x$ in $[0,1]$
such that $\mu_x(\e_1=0)=\a$. \newline
2. A typical $x$ from the viewpoint of Lebesgue measure, of
course, corresponds to the \linebreak case {\bf a.} of the theorem.
\endremark
\subhead Examples \endsubhead
We illustrate the possible situations in the previous
theorem with four examples. For a better illustration we will use
the following convention:

\bigskip
\epsfysize=3cm
\centerline{\epsfbox{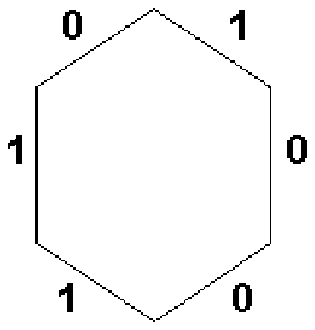}}
\bigskip\bigskip
\nopagebreak

\centerline{Fig. 4}
\medskip

We will use generalized canonical expansions, writing
$x\sim(\e_1 \e_2 \dots)$.\newline
{\bf 1.} $x=\la^{-2}\sim (1000\dots)$. Here $\phi_x$ is countable and isomorphic
to the stationary Markov compactum with the matrix $\(\smallmatrix 1&1\\0&1
\endsmallmatrix\)$. This compactum consists of the sequence $0^\infty$
and the sequences $0^k1^\infty$ for $k\ge0$. For a central measure
$\mu_x$, $\mu_x(0^k1^\infty)=\mu_x(0^l1^\infty)$ for any $k,l$, hence,
this measure is concerntrated on the path $0^\infty$ corresponding to
the path $(010101\dots)$ in the initial compactum (see Fig.~5 below).

\bigskip
\epsfysize=5cm
\centerline{\epsfbox{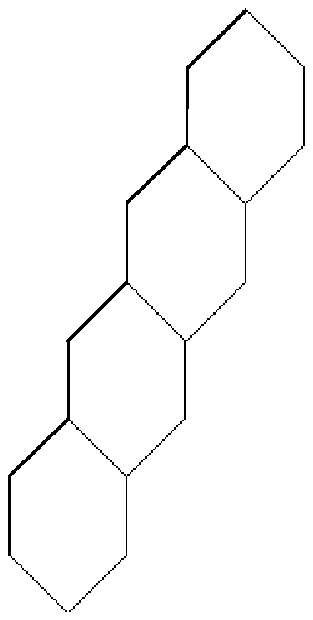}}
\bigskip\bigskip
\nopagebreak

\centerline{Fig. 5. The case $x=\la^{-2}$}
\medskip

{\bf 2.} $x=\frac12\sim(100)^\infty$. Here $B^{(0)}=\emptyset,\,
B^{(1)}=\(1(00)^1\)^\infty$. We have $\phi_x=\prod_1^\infty\{011, 100\}$,
and thus $T_x$ is isomorphic to the 2-adic shift, i.e. the shift by 1
in the group of dyadic integers. Therefore, $T_x$ has the binary rational
purely discrete spectrum. We depict the way of recoding the paths
in $\phi_x$ into the full dyadic compactum by the rule ``$011\sim 0,\,\,
100\sim 1$''(see Fig.~6).

\bigskip
\epsfysize=5cm
\centerline{\epsfbox{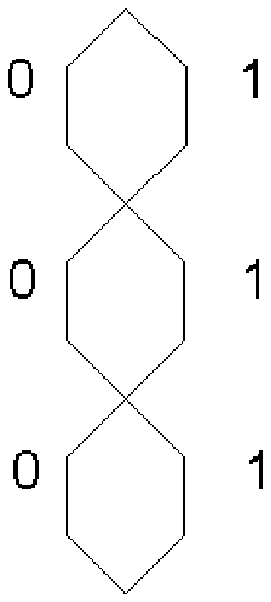}}
\bigskip\bigskip

\nopagebreak
\centerline{Fig. 6. Recoding the paths for $x=\frac12$, case {\bf a}}
\medskip

{\bf 3.} $x\sim (1(0001)^\infty)$. Here $\a=[1,1,1,\dots]=\la^{-1}$, and
$T_x$ acts as the rotation by the golden ratio. Fig.~8
shows the way of recoding the paths in $\phi_x$ into the usual model
for this rotation (``Fibonacci compactum''). Note that the natural ordering
of these paths is alternating (see Fig.~7, 8).

\bigskip
\epsfysize=6cm
\centerline{\epsfbox{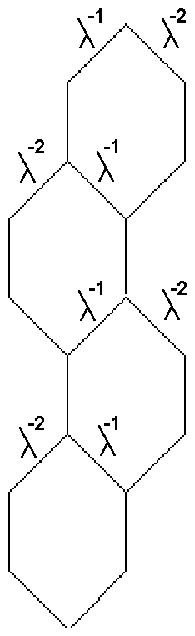}}
\bigskip\bigskip
\nopagebreak

\centerline{Fig. 7. Transition measures for $\a=\la^{-1}$}
\medskip

\bigskip
\epsfysize=6cm
\centerline{\epsfbox{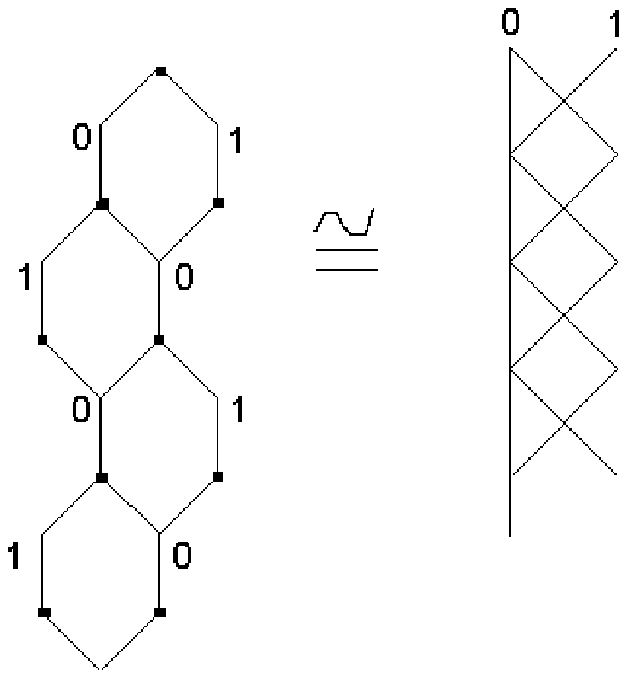}}
\bigskip
\nopagebreak

\centerline{Fig. 8. Recoding the paths for $\a=\la^{-1}$, case {\bf b}}
\medskip

{\bf 4.} $x\sim \(1001(0001)^\infty\)$. For this $x$, the transformation
$T_x$ acts as the special automorphism over the rotation by $\la^{-1}$
with a single step equal to the base (see Fig.~9).

\bigskip
\epsfysize=6cm
\centerline{\epsfbox{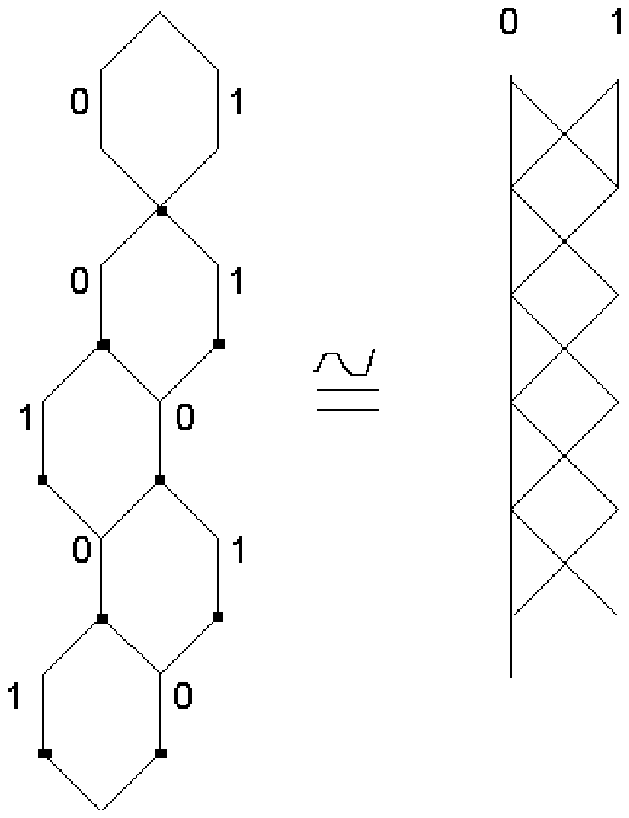}}
\bigskip

\nopagebreak
\centerline{Fig. 9. Recoding the paths for case {\bf c}}
\medskip

\head Appendix B. Arithmetic expression for block expansions \endhead

We recall that in Section~2 we have defined the mapping $\Psi$ assigning
to a regular $x\in(\la^{-1},1)$ the sequence of blocks
$B_1(x), B_2(x),\dots$. In this Appendix
we are going to specify
the mapping $\Psi$ in an arithmetic way. To this end, we gather
the canonical expansion of a given regular $x$ blockwise.

Recall that similarly to the canonical expansion~(1.1) of reals,
there exists the corresponding representation of positive integers.
Namely, each $N\in\Bbb N$ has a unique representation in the form
$$
N=\sum_{i=1}^k \e_i F_i,
$$
where $\e_i\in\{0,1\},\, \e_i\e_{i+1}=0,\, \e_k=1$
for some $k\in\Bbb N$. It is usually
called the {\it Zeckendorf} decomposition. We denote by $\Cal F$
the class of positive integers whose Zeckendorf decomposition has
$\e_1=1$ and $\e_i\equiv 0$ for all even $i$. Obviously, $\Cal F$ as
a subset of $\Bbb N$ has zero density. Let the {\it height} of
$e$ with a finite canonical expansion of the form
$e=\sum_j \e_j\la^{-j}$ be, by definition,
the positive integer $h(e):=\max\{j:\e_j=1\}$.
\proclaim{Proposition B.1} Each regular $x\in(\la^{-1},1)$ has
a unique representation of the form
$$
x=\sum_{j=1}^\infty e_j(x)\la^{-\sum\limits_{i=1}^j n_i(x)}, \tag B.1
$$
where
\roster
\item "(i)" $e_j(x)=m\la-n,\, n\in \Cal F,\, n=[m\la]$.
\item "(ii)" $n_j$ is odd, $n_j\ge 3$ for $j\ge 1$.
\item "(iii)" $n_j\ge 2h(e_j)+1$ for all $j$.
\endroster
\endproclaim
\demo{Proof} Let $\Psi(x)=B_1 B_2\dots$ be the block expansion
of $x$.
Suppose $B_j=B_j(a_1^{(j)},\dots,a_{t_j}^{(j)})$. We set
$n_j(x):=2\sum_{i=1}^{t_j}a_i^{(j)}+1,\,\, j\ge 1$,
i.e. $n_j$ is the length of the $j$'th block. Let $m_0^{(j)}=0,\,
m_k^{(j)}=\sum_{i=1}^ka_i^{(j)}$, and let
$e_j$ be the ``value'' of $B_j$ in the sense of formula~(1.1)
as if it were the first block, i.e.
$$
\align
e_j(x)&=\la^{-1}+\la^{-1}\sum_{k=1}^{\frac12(t_j-1)}
\sum_{\nu=1}^{a_{2k}^{(j)}}\la^{-2\(m_{2k-1}^{(j)}+\nu\)} \\
      &=\la^{-1}+\la^{-2}\sum_{k=1}^{\frac12(t_j-1)}
\(\la^{-2m_{2k-1}^{(j)}}-\la^{-2m_{2k}^{(j)}}\), \quad t_j\,\, \text{odd}, \\
e_j(x)&=\la^{-1}+\la^{-1}\sum_{k=1}^{\frac12 t_j}
\sum_{\nu=1}^{a_{2k+1}^{(j)}}\la^{-2\(m_{2k-2}^{(j)}+\nu\)} \\
      &=\la^{-1}+\la^{-2}\sum_{k=1}^{\frac12 t_j}
\(\la^{-2m_{2k-2}^{(j)}}-\la^{-2m_{2k-1}^{(j)}}\), \quad t_j\,\, \text{even}.
\endalign
$$
The uniqueness of expansion~(B.1) follows from the condition~(iii)
and from the uniqueness of expansion~(1.1) for any finite sequence
(and, therefore, for any block). The fact that $n\in\Cal F$ follows
from the definition of block.
\enddemo
\definition{Definition} We call the expansion of $x\in (\la^{-1},1)$
of the form~(B.1) satisfying the conditions (i)--(iii) the {\it arithmetic
block expansion}.
\enddefinition

\remark{Remark {\rm 1}}  $n_j$ and $e_j$ depend on $B_j$ only.
\endremark
\remark{Remark {\rm 2}}
In fact, series~(B.1) is nothing but series~(1.1) rewritten in a different
notation. However, we will see that it has its own dynamical sense
(see relation~(B.2) below).
\endremark
\remark{Remark {\rm 3}} By Item~(iii), the quantities $e_j$ and
$n_j$ are not completely independent. Let $l_j:=h(e_j)$.
Taking into consideration new quantities $s_j:=n_j-2l_j$ and
representing $n_j$ as the sum $s_j$ and $2l_j$ in formula~(B.1),
we come to independent multipliers, but this new form of the block
representation does not seem to be natural.
\endremark
\remark{Remark {\rm 4}} In terms of
arithmetic block expansions the goldenshift acts as
$$
\S(x)=\sum_{j=2}^\infty e_j(x)\la^{-\sum\limits_{i=2}^jn_i(x)}.\tag B.2
$$
\endremark

\head Appendix C. Computation of densities and the polymorphism $\Pi$
\endhead

We return to the subject of the first section. Recall that we have
already denoted
the transformation $x\mapsto \{\la x\}$ by $T$, and $R$ stands
for the rotation of the circle $\Bbb R/\Bbb Z$ by the angle $\la^{-1}$.

\subhead C.1. Computation of densities \endsubhead

\proclaim{Proposition C.1} The densities $\frac{d(R\mu)}{d\mu},\,
\frac{d(\tau\mu)}{d\mu}$ and $\frac{d\nu}{d\mu}$ are unbounded
and piecewise constant with a countable number of steps.
\endproclaim
\demo{Proof} By virtue of the results of Section~1, it suffices to prove
the proposition only for $\frac{d(R\mu)}{d\mu}$.
Let $E$ be a Borel subset of $(0,1)$. The idea of the study
lies in the fact that $R$ does not change any block beginning with
the second. As usual, we consider three cases. \newline
{\bf I.} $E\subset (0, \la^{-2})$. If $E\subset (\la^{-2k},\la^{-2k+1}),
\, k\ge 1$, then each point $x$ of the set $E$ has the canonical
expansion (1.1) of the form $0^{2k-1}10*$. Hence the canonical expansion
of $x+\la^{-1}$ is $1(00)^{k-1}10*$, and
$$
\frac{\mu(E+\la^{-1})}{\mu E}\equiv k,
$$
as $f(1(00)^{k-1})=k$. If, on the contrary,
$E\subset (\la^{-2k-1},\la^{-2k}),\ k\ge1$, then the situation is as follows.
This interval in terms of the canonical expansion is
$\bigcup_{\ov B} 0^{2k} \ov B\mod 0$,
where the union runs over all closed blocks $\ov B$.\footnote{We say that
a 0-1 word is a {\it closed block} if it has the form
$B1$ for some block $B$.}
We have two subcases. \newline
{\bf Ia.} Let in terms of the canonical expansion,
$E\subset 0^{2k}1(00)^{a_1}(01)^{a_2}\dots(00)^{a_t}1$. Here
$E+\la^{-1}\subset 1(00)^{k-1}(01)(00)^{a_1}(01)^{a_2}\dots(00)^{a_t}1$,
hence $\frac{\mu(E+\la^{-1})}{\mu E}=\frac{f(B')}{f(B)}$, where
$B'$ is the closed block defined as
$B'=B'(k-1,1,a_1, a_2,\dots, a_t)$. So, we conclude from Lemma~2.1 that
$$
\frac{\mu(E+\la^{-1})}{\mu E}=\frac{kp+(k+1)q}{p+q},
$$
where, as usual, $\frac pq=[a_1, a_2,\dots, a_t]$. \newline
{\bf Ib.} In the same terms, suppose
$E\subset 0^{2k}1(01)^{a_1}(00)^{a_2}\dots(00)^{a_t}1$. Similarly
to the above,
$$
\frac{\mu(E+\la^{-1})}{\mu E}=\frac{(k+1)p+kq}{p+q}.
$$
{\bf II.} Let
$E\subset (\la^{-2},\la^{-1})$. This case is analogous to Case~I.
If $E\subset (\la^{-2}+\la^{-2k-3},\la^{-2}+\la^{-2k-2}),\, k\ge 1$, then
$$
\frac{\mu(E-\la^{-2})}{\mu E}\equiv \frac 1k.
$$
If $E\subset (\la^{-2}+\la^{-2k-2},\la^{-2}+\la^{-2k-1}),\, k\ge 1$, then
$$
\frac{\mu(E-\la^{-2})}{\mu E}\equiv
\cases
\frac{p+q}{kp+(k+1)q},
            & E\subset 01(00)^{k-1}01(00)^{a_1}(01)^{a_2}\dots(00)^{a_t}1\\
\frac{p+q}{(k+1)p+kq},
            & E\subset 01(00)^{k-1}01(01)^{a_1}(00)^{a_2}\dots(00)^{a_t}1.
\endcases
$$
{\bf III.} Let $E\subset (\la^{-1},1)$.  If $E\subset(\la^{-1},\la^{-1}+
\la^{-4})$, then $\mu(E-\la^{-2})=\mu E$.
If $E\subset(\la^{-1}+\la^{-4},1)$, then
$E-\la^{-2}\subset(\la^{-2},\la^{-1})$, hence $E-\la^{-2}\subset
010*$.   \newline
{\bf IIIa.} Let
$E-\la^{-2}\subset 1(00)^{a_1}(01)^{a_2}\dots(00)^{a_t}1$, then
$$
\frac{\mu(E-\la^{-2})}{\mu E}=1+\frac pq.
$$
{\bf IIIb.} Let
$E-\la^{-2}\subset 1(01)^{a_1}(00)^{a_2}\dots(00)^{a_t}1$. Here
$$
\frac{\mu(E-\la^{-2})}{\mu E}=1+\frac qp.
$$
The proof is complete.
\enddemo
\remark{Remark {\rm1}} Let
$d=\frac{d(R\mu)}{d\mu}$. Then by relation (1.7),
$\frac{d(\tau\mu)}{d\mu}=\frac12(d+1)$, and by formula~(1.12),
$$
\frac{d\nu}{d\mu}(x)=
\cases
\frac23+\frac13d(x)+\frac16 d^{-1}(x+\la^{-1}),
             & x\in [0, \la^{-2}) \\
\frac23+\frac13d(x),
             &x \in [\la^{-2}, \la^{-1}) \\
\frac12 +\frac13d(x),
             &x \in [\la^{-1}, 1].
       \endcases
$$
\endremark
\remark{Remark {\rm2}} From this relation follows Proposition~2.10.
\endremark

\subhead C.2. The polymorphism $\Pi$ \endsubhead
Let as above $\sigma:\Sigma\to\Sigma$ be the one-sided shift. Let us ask the
natural question: what is the image of $\sigma$ on the interval
$[0,1]$ under the mapping~$\pi$ defined by the formula~(3.1)?

Note first that the partition into $\pi$-preimages of singletons
is not invariant under $\sigma$. Indeed, if, say, $x=0110^\infty$, then
$\sigma\n(x)=0^\infty$, while $\n\sigma(x)=110^\infty$. Thus,
$\Pi:=\pi\sigma\pi^{-1}:[0,1]\to[0,1]$ is
a multivalued mapping,
i.e. a {\it polymorphism} by the terminology of \cite{Ver4}.

Recall that a {\it measure-preserving polymorphism}
of a measure space $(X,\goth A,\mu)$ is the diagram
$$
(X,\mu) @<{\pi_1}<<(Y,\nu) @>{\pi_2}>>(X,\mu),
$$
where $\pi_1,\pi_2$ are homomorphisms of measure spaces such that
$\pi_i\nu=\mu,\ i=1,2$. Instead of an arbitrary $Y$
it suffices to consider $X\times X$ with the coordinate projections
and a certain ``bistochastic" measure $\nu$ (i.e. a measure
on the sigma-algebra $\goth A\times\goth A$ with given marginal
measures). Such a polymorphism is called {\it reduced}.

If $T$ is an automorphism of the space $Y$ with an invariant
measure $\ga$ and $\zeta$ is a measurable partition, one can
define the polymorphism $T_\zeta:(Y_\zeta,\ga_\zeta)$ into itself
as follows. Consider two partitions of $Y$, namely,
$\zeta$ and $T^{-1}\zeta$. Identifying $Y_\zeta$ and $Y_{T^{-1}\zeta}$
in the natural way, we obtain the diagram
$$
(Y_\zeta,\nu_\zeta)\longleftarrow(Y,\nu)\longrightarrow(Y_\zeta,\nu_\zeta).
$$
Let $\xi=\zeta\vee T^{-1}\zeta$. Then $Y_\xi\subset Y_\zeta\times
Y_{T^{-1}\zeta}=Y_\zeta\times Y_\zeta$, and the reduced automorphism
$$
(Y_\zeta,\nu_\zeta)\longleftarrow(Y_\xi,\nu_\xi)\longrightarrow
(Y_\zeta,\nu_\zeta)
$$
can be easily interpreted: it is a Markov (multivalued) mapping of the
factor space $Y_\zeta$ with the invariant measure $\nu_\zeta$.
If $\zeta$ is a $T$-invariant partition (i.e. $T^{-1}\zeta\prec\zeta$),
then the polymorphism is the factor endomorphism
$T_\zeta:Y_\zeta\to Y_\zeta$. That is why the polymorphism in this
context is a generalization of the notion of endomorphism (see
\cite{Ver4}).

We are going to make use of these notions in our case.
We define the polymorphism $\Pi$ as
the subset of $[0,1]\times[0,1]$ defined as
$$
\Pi(x)=\cases \la x,                   & 0\le x<\la^{-2}\\
              \la x\cup\la x-\la^{-1}, & \la^{-2}\le x<\la^{-1}\\
              \la x-\la^{-1},          & \la^{-1}\le x\le1
       \endcases
$$
and provided with the measure $\ov\mu$ which is the image
of the product measure $p$. By definition
of the polymorphism and the Erd\"os measure,
$\Pi\mu=\Pi^{-1}\mu=\mu$, and $\mu$
being the projection of $\ov\mu$ to both axes.

\bigskip
\epsfysize=5cm
\centerline{\epsfbox{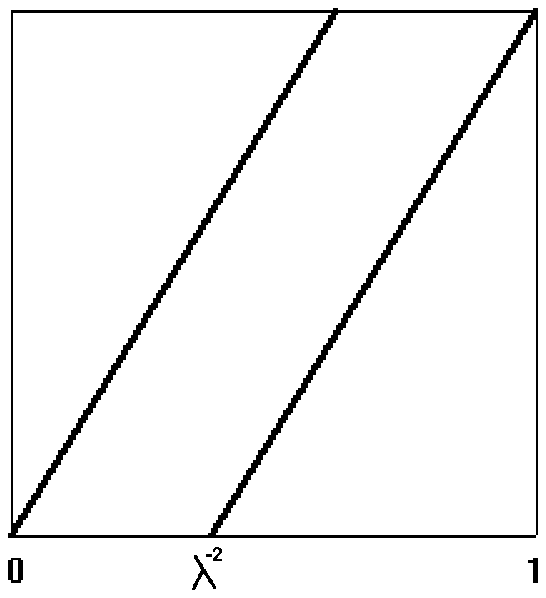}}
\bigskip
\nopagebreak

\centerline{Fig. 10. The polymorphism $\Pi$}
\medskip

This polymorphism was considered in \cite{VerSi}.
Note that for a Borel set $E$,
$\Pi^{-1}E=\pi\sigma^{-1}\pi^{-1}E=
\la^{-1}E\cup\la^{-1}E+\la^{-2}$.
Let $\ga=(\ga_1,\ga_2)$ be the corresponding partition of $\Pi$.
It is possible to show that there exists a countable partition
of $[0,1]$ into the intervals $\{G_k\}_{k=1}^\infty$
such that for any $k$ and any $G\subset G_k$ the ratio
$\ov\mu(\Pi^{-1}G\cap\ga_1)/\ov\mu((\Pi^{-1}G\cap\ga_2)$ is constant.
In particular, if $G\subset(\la^{-2n},\la^{-2n+1})$ for $n\ge1$,
then this ratio is equal to $n$ (specifically, for
$G\subset(\la^{-1},\la^{-2})$, it equals 1). The method of the
proof is the same as in Proposition~C.1.

\head Appendix D. An independent proof of Alexander-Zagier's formula \endhead

In this appendix we will present the second proof of formula~(3.3) which
reveals some new relations between certain structures of the Fibonacci
graph~$\Phi$.

We first recall that the quantity $f_n(k)$ is nothing but the frequency
of the $k$'th vertex on the $n$'th level of the Fibonacci graph
which was denoted by $D_n$
(see the beginning of Section~3). We have $\#D_n=F_{n+2}-1$.

Consider level $n$ of the Fibonacci graph for $n=2N+1$.
We denote the {\it middle}
part of $D_n$, i.e. the segment from $F_n$ to $F_{n+1}-1$, by $D'_n$.
The Erd\"os measure of $D'_n$
obviously equals $\frac13+O(\la^{-n})$, and we will introduce the
partition of $D'_n$ into $2^{N-1}-1$ intervals of vertices in
the following way.

Recall that a  Euclidean vertex
is one, where the first block can end (they are
marked in Fig.~3 for $D_3$ and $D_5$) and that these vertices form
the Euclidean binary tree (see Section~3).
There are $2^{N-1}$ such vertices at level~$2N+1$, and all of them
lie in $D'_n$. Let $V_k^{(N)}$
denote the $k$'th Euclidean vertex from the left on level $2N+1$.
\definition{Definition} An open interval of vertices
$\Omega_k^{(N)}:=(V_k^{(N)}, V_{k+1}^{(N)})$ will
be called a {\it Euclidean} interval.
\enddefinition

So, we have divided the set of vertices $D'_n$ into $2^{N-1}$
Euclidean vertices $\(V_k^{(N)}\)_{k=1}^{2^{N-1}}$
and $2^{N-1}-1$ open  intervals $\Omega_k^{(N)}$. Now
we introduce the subgraph $\Gamma_V$ associated with each Euclidean
vertex $V$. It is defined as the one containing all the successors
of $V$ in the sense of the Fibonacci graph, except any other
Euclidean vertices (see Fig.~11).

\bigskip
\epsfysize=6cm
\centerline{\epsfbox{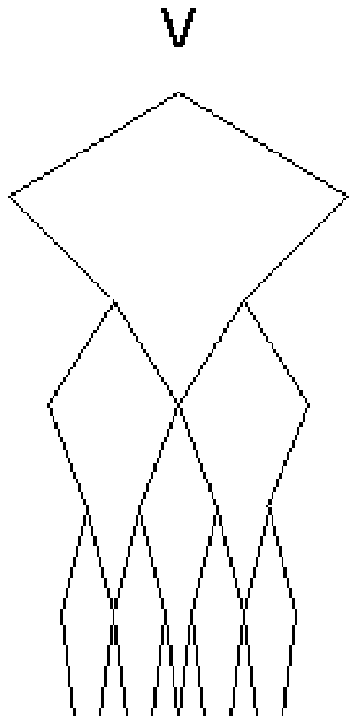}}
\bigskip
\nopagebreak
\centerline{Fig. 11. The graph $\Gamma_V$}
\medskip

We state a straightforward lemma.
\proclaim{Lemma D.1} For any Euclidean interval $\Omega_k^{(N)}$
there is a unique Euclidean vertex $V_i^{(j)},\,j<N$, such that
$\Omega_k^{(N)}\subset \Gamma_{V_i^{(j)}}$.
\endproclaim
So, any Euclidean interval is determined by a certain Euclidean vertex
on one of the preceding odd levels of the Fibonacci graph. Moreover,
in the notation of the above lemma,
the entropy of $\Omega_k^{(N)}$ may be computed in terms of the frequency
of $V_i^{(j)}$ and the entropy of $\wt{D}_{2N+1-j}$. Namely,
let $H_n:=\sum\limits_{k=F_n}^{F_{n+1}-1} f_n(k)\log_\la f_n(k)$,
and let next
$H_{2N+1}=\sum_{j=1}^N H_{2N+1}^{(j)}$, where $H_{2N+1}^{(j)}$
denotes the sum over the vertices $V\in\Gamma_{V_i^{(j)}}$ for
all $i\le 2^{j-1}$. So, $H_{2N+1}^{(j)}$ corresponds to all
Euclidean vertices of level~$2j+1$.

Let next $\varphi_i^{(j)}$
denote the frequency of $V_i^{(j)}$. For instance, for $j=3$,
$\, \varphi_1^{(3)}=\varphi_4^{(3)}=4,\,\, \varphi_2^{(3)}=
\varphi_3^{(3)}=5$. In this notation $k_j=\sum_{i=1}^{2^{j-1}}
\varphi_i^{(j)}\log_\la\varphi_i^{(j)}$.

For any $v\in (V_k^{(N)},V_{k+1}^{(N)})$ its frequency equals $f(V_i^{(j)})$
times the frequency of the corresponding vertex of the central part of
level~$2N+1-j$. So, we established an essential relationship between
the central part of level~$2N+1$ of the Fibonacci graph and all
Euclidean vertices $V_i^{(j)},\,1\le j\le N,\,1\le i\le 2^{j-1}$.
\proclaim{Lemma D.2} The following recurrence relation holds:
$$
H_{2N+1}=\frac23\sum_{j=1}^{N-1} 3^j H_{2N-2j}+\frac13\cdot 4^N\cdot
\sum_{j=1}^N \frac{k_j}{4^j} + O\(\sum_{j=1}^N k_j\),
\quad N\to\infty. \tag D.1
$$
\endproclaim
\demo{Proof} By the above considerations,
$$
\align
H_{2N+1}^{(j)}&=\sum_{i=1}^{2^{j-1}}\sum_{k=F_{2N-2j}}^{F_{2N-2j+1}-1}
\varphi_i^{(j)}f_{2N-2j}(k)\log_\la\(\varphi_i^{(j)}f_{2N-2j}(k)\) \\
              &=\sum_{i=1}^{2^{j-1}}\varphi_i^{(j)}
\(H_{2N-2j}+\frac13\log_\la\varphi_i^{(j)}\cdot\(4^{N-j}+O(1)\)\) \\
              &=2H_{2N-2j}\cdot 3^{j-1}+\frac13\cdot\frac{k_j}{4^j}\cdot 4^N
+O(k_j)
\endalign
$$
(we used the fact that  $\sum_{i=1}^{2^{j-1}}\varphi_i^{(j)}=2\cdot 3^{j-1}$
easily obtained from Proposition~2.9). Hence relation~(D.1) follows.
\enddemo
\remark{Remark} Formula (D.1) shows that the entropy of the $n$'th
level with $n$ odd can be computed by means of the entropies of
the previous even levels and the entropy of the Euclidean
tree.
\endremark
Now we are ready to complete the second proof of formula~(3.3). We have
$$
nH_\mu \sim \sum_{k=0}^{F_{n+2}-2}\frac{f_n(k)}{2^n}
\log_\la\frac{2^n}{f_n(k)},
$$
whence
$$
nH_\mu \sim 3\sum_{k=F_n}^{F_{n+1}-1}\frac{f_n(k)}{2^n}
\log_\la\frac{2^n}{f_n(k)},
$$
and
$$
H_n \sim \frac13\bigl(\La-H_\mu\bigr)n2^n. \tag D.2
$$
From relation~(D.2) it follows that in the sum
$\sum\limits_{j=1}^{N-1} 3^jH_{2N-2j}$ the first terms are more valuable than
the last. Thus, from formulas~(D.1) and (D.2) and from
the fact that $k_N=\sum_{i=1}^{2^{N-1}}\varphi_i^{(N)}\log_\la\varphi_i^{(N)}
<2(N-1)3^{N-1}$ it follows that
$$
\align
\frac13\bigl(\La-H_\mu\bigr)\cdot(2N+1)2^{2N+1}&\sim
\frac23\sum_{j=1}^N 3^j\cdot \frac13\bigl(\La-H_\mu\bigr)
(2N-2j)4^{N-j} \\
&+\frac13\cdot 4^N\sum_{j=1}^N\frac{k_j}{4^j},
\endalign
$$
whence, after straightforward computations,
$$
18\bigl(\La-H_\mu\bigr) \sim \sum_{j=1}^N\frac{k_j}{4^j},\quad
N\to\infty. \qed
$$
\remark{Remark} The Euclidean tree, being symmetric,
naturally splits into two binary
subtrees (left and right) being symmetric. If we label each vertex of
the left subtree with the corresponding rational $p/q$, then this left
subtree turns out to coincide with the {\it Farey} tree introduced and
studied in detail in \cite{Lag}.
\endremark

\enddocument